\documentclass[10pt,twoside]{amsart}
\usepackage{amsmath, amsfonts, amsthm, amssymb} 
\newtheorem{definition}{Definition}[section]
\newtheorem{lemma}[definition]{Lemma}
\newtheorem{theorem}[definition]{Theorem}
\newtheorem{proposition}[definition]{Proposition}

\newtheorem{remark}[definition]{Remark}

\numberwithin{equation}{section} 
\newcommand \be   {\begin{equation}}
\newcommand \ee   {\end{equation}} 
\newcommand \RR      	{\mathbb{R}} 
\newcommand \eps      \epsilon  
\newcommand \ubar   	{{\overline u}}

\newcommand \del     	\partial 
\newcommand \MM     	{M}  
\newcommand \Sbold   {S} 
\newcommand \nn       	{\mathbf n}

\newcommand \ttt  	{\mathbf t} 
\newcommand \Lcal  	{\mathcal L} 
\newcommand \Ocal  	{\mathcal O} 
  
\newcommand \TS		{T_x S^2}

\newcommand \TM {T_x M}
\newcommand \coTM {T_x^\star M}

\newcommand \nablac  {\nabla_g \! \cdot \!} 
\newcommand \restrict	[1]{\raisebox{-1mm}{\ensuremath{|_{#1}}}}
\newcommand \Lie 	{\mathcal{L}}
  
\newcommand \TT		{\mathcal{T}}

\newcommand{\dK}{\partial K}
\newcommand{\nek}{\nn_{e,K}}

\newcommand{\fek}{f_{e,K}}
\newcommand{\feke}{f_{e,K_e}}
\newcommand{\Fek}{F_{e,K}}
\newcommand{\Feke}{F_{e,K_e}}
\newcommand{\Lip}{\mathop{\mathrm{Lip}}}
\newcommand{\unk}{u^n_K}
\newcommand{\unke}{u^n_{K_e}}
\newcommand{\unnk}{u^{n+1}_K}
\newcommand{\unnkee}{u^{n+1}_{K,e}}
\newcommand{\unnkeee}{u^{n+1}_{K_e,e}}

\newcommand{\tunnkee}{\tilde u^{n+1}_{K,e}}
\newcommand{\rnnkee}{R^{n+1}_{K,e}}
\newcommand{\sumke} {\sum_{\substack{K \in \TT^h \\ e \in \dK}}}
\newcommand{\sumk} {\sum_{K\in\TT^h}}
\newcommand{\pne}{\phi_e^n}
\newcommand{\phnk}{\hat\phi^n_K}
\newcommand{\pnk}{\phi_K^n}
\newcommand{\inttn}{\int_{t_n}^{t_{n+1}}}

\begin{document} 
\bibliographystyle{plain}  
\title[\bf Hyperbolic Conservation Laws on Manifolds]
{Hyperbolic Conservation Laws on Manifolds. 
\\
Total Variation Estimates and 
\\
the Finite Volume Method}
\author[Amorim, Ben-Artzi, LeFloch]
{
Paulo Amorim$^1$, Matania Ben-Artzi$^2$
\\
\and
\\ 
Philippe G. LeFloch$^1$
}  
\thanks{Completed on June 2005. To appear in Methods and Applications of Analysis. 
\newline
$^1$ Laboratoire Jacques-Louis Lions 
\& Centre National de la Recherche Scientifique, UMR 7598, 
Ê      University of Paris 6, BC 187, 75252 Paris, 
 Ê     France. E-mail : {\tt amorim@ann.jussieu.fr, lefloch@ann.jussieu.fr.} 
 \newline 
$^2$ Institute of Mathematics, Hebrew University, Jerusalem 91904, Israel. E-mail: {\tt mbartzi@math.huji.ac.il}
\newline
2000\textit{\ AMS Subject Classification:} 35L65, 74J40, 58J, 76N10. 
}  

\begin{abstract} This paper investigates some properties of entropy solutions of  hyperbolic conservation laws 
on a Riemannian manifold. First, we generalize the Total Variation Diminishing (TVD) property to manifolds, 
by deriving conditions on the flux of the conservation law and a given vector field 
ensuring that the total variation of the solution along the integral curves of the vector field 
is non-increasing in time. Our results are next specialized to the important case of a flow on the $2$-sphere, 
and examples of flux are discussed. Second, we establish the convergence of the finite volume methods
based on numerical flux-functions satisfying monotonicity properties. Our proof requires detailed estimates 
on the entropy dissipation, and extends to general manifolds
an earlier proof by Cockburn, Coquel, and LeFloch in the Euclidian case. 
\end{abstract}

\maketitle

\section{Introduction} 
\label{IN-0}

In this paper, following \cite{BenArtziLeFloch}, we investigate some properties of entropy solutions 
to the Cauchy problem for hyperbolic conservation laws on a manifold~: 
\be 
\del_t u + \nablac f(u, \cdot) = 0,  
\qquad u=u(t,x) \in \RR, \quad t \geq 0, \, x \in M. 
\label{IN.1}
\ee
Here, $M$ is a $d$-dimensional, smooth manifold endowed with a Riemannian metric $g$ 
and, for each constant $\ubar \in \RR$, the map $x \mapsto f(\ubar,x)$ is a smooth vector field on $M$, 
i.e.~a section of the tangent bundle $TM$. 
The well-posedness theory for the hyperbolic conservation law was recently 
established in Ben-Artzi and LeFloch \cite{BenArtziLeFloch}. 
As in Kruzkov's theory \cite{Kruzkov} which applies to conservation laws in the Euclidian setting $M = \RR^d$, 
one is interested in weak solutions of \eqref{IN.1} in the sense of distributions that 
are constrained by \emph{entropy inequalities.} For instance, in the case that $f$ is divergence-free, that is,
\be 
\nablac f(\ubar, \cdot) = 0 
\quad \text{ for every constant } \ubar, 
\label{IN.2}
\ee 
the entropy inequalities read
\be 
\del_t U(u) + \nablac F(u, \cdot) \leq 0,
\label{IN.3}
\ee 
for every convex $U: \RR \to \RR$, where the vector field $x \mapsto F(\ubar,x)$ is 
$$
\del_u F(\ubar,x) := \del_u U(\ubar) \, \del_u f(\ubar,x), \quad \ubar \in \RR, \, x \in M.  
$$

Let us emphasize that 
the equation \eqref{IN.1} is a geometric partial differential equation which naturally 
depends on the geometry of the manifold, only. In particular, all estimates derived on solutions to
\eqref{IN.1} should take a form that is completely independent of any particular system of local coordinates of $M$,
although in practice, for the proofs it will often be convenient to introduce a particular chart to represent the manifold.
Throughout, the convention of implicit summation over repeated indices will be used. 
In local coordinates $x=(x^j)_{1 \leq j \leq d}$, we will use the short-hand notation $\del_j := \del/\del x^j$.  
Recall that the divergence operator arising in \eqref{IN.1} takes the form 
$$
\nablac f(u(t,x), x)  := {1 \over \sqrt{|g(x)|}} \, \del_j \big(\sqrt{|g(x)|} \, f^j(u(t,x),x) \big), 
$$
where  $(g_{ij})$ are the coordinates of the metric tensor $g$ and $|g| := \det (g_{ij})$. Here, $f^j(\ubar,x)$ are, 
for every $\ubar$, the coordinates of the vector field $f(\ubar, \cdot)$ at the point $x$. 

The present paper supplements the well-posedness results established in the companion paper \cite{BenArtziLeFloch}
and has two main objectives. 

First, we investigate some properties of the total variation of solutions of \eqref{IN.1} on a manifold. 
Precisely, we generalize the Total Variation Diminishing (TVD) property to manifolds, by 
deriving conditions on the flux $f$ and a given vector field $X$ 
ensuring that the total variation of the solution along the integral curves of $X$ 
is non-increasing in time.  
Recall that TVD schemes for nonlinear hyperbolic problems 
play a central role in scientific computation, for instance in gas dynamics. 
Our diminishing total variation properties provide certain a~priori estimates which, for instance, can be tested numerically
and may help in designing robust schemes that are consistent with large-time asymptotics.   
Note that solutions of a conservation law on a manifold need not have a total variation that is bounded uniformly in time.
The geometric effets may contribute to amplify the wave strengths, and 
the total variation of a solution may blow-up in the large as $t \to \infty$.  
We also investigate here in some detail the case where the manifold is the sphere $\Sbold^2$
embedded in $\RR^3$.
Examples of flux (as vector fields in $\RR^3$ tangent to the sphere)
are discussed from both the intrinsic and embedded standpoints. 

Second, we consider the numerical approximation of the entropy solutions of \eqref{IN.1} 
and we establish the convergence of the finite volume scheme 
on manifolds when the numerical flux-functions are depend monotonically upon their arguments. 
Our result is an extension to general manifolds of a theorem due to  
Cockburn, Coquel, and LeFloch \cite{CCL} (see also \cite{CCLS}) in the Euclidian case. While the convergence 
of the finite volume scheme is easily established on a cartesian mesh, due to the invariance by translation of 
both the equation and the mesh, the convergence proof for non-cartesian meshes is more involved
and, as was pointed out in \cite{CCL}, requires DiPerna's concept of measure-valued solutions \cite{DiPerna}. 
To handle \eqref{IN.1}, we therefore must rely on Ben-Artzi--LeFloch's extension to manifolds \cite{BenArtziLeFloch} 
of DiPerna's theorem.

An outline of this paper is as follows. In Section~2, we 
derive conditions on the flux and the vector field $X$ ensuring that 
the function $X(u)$ satisfies a conservation law whenever $u$ is a solution of \eqref{IN.1}. 
This leads us to a diminishing total variation property for a class of flux and vector fields. 
In Section~3, we discuss the structure and general properties of conservation laws on the sphere $S^2$.  
In Section~4, we apply the framework of Section~2 to the case of the sphere and derive total variation bounds.
Next, Sections~5 and 6 are devoted to the statement and to the proof of 
the convergence of the finite volume scheme. 

\section{Total variation diminishing estimates on a general manifold}   
\label{BV-0}

In the present section, we derive conditions on the flux $f$ of the conservation law \eqref{IN.1}
which ensure that the total variation (at least along certain vector fields)
of the entropy solutions of \eqref{IN.1} is non-increasing in time. We will first state the main results
(Proposition~\ref{BV-1} and Theorem~\ref{BV-3}) 
and then recall some elementary notions  from differential geometry, before giving the proofs
of the results.

\subsection{Statement of the estimate} 
Throughout, $M$ is a smooth, $d$-dimensional, Riemannian manifold $(M,g)$, which has
no boundary and need not be compact. The following result provides us with a key identity on 
directional derivatives of solutions of \eqref{IN.1}.

\begin{proposition} 
\label{BV-1} 
Let $u: \RR_+ \times M \to \RR$ be a smooth solution of the conservation law \eqref{IN.1} on $M$ where,  
$f=f(\ubar, \cdot)$ is a (smooth) vector field depending on the parameter $\ubar$. 
Then,  given any (smooth) vector field $X$ on $M$ the function $w := X(u): \RR_+ \times M \to \RR$ 
satisfies the (linear) hyperbolic equation 
\be 
\del_t w + \nablac \big( w\,  f_u(u, \cdot) \big) = - g \big(\nabla_g u, (\Lie_X f_u)(u, \cdot) \big) 
       - X(\nablac f)(u,\cdot), 
\label{BV.1} 
\ee 
where $(\Lie_X f_u)(u, \cdot)$ denotes the Lie bracket of the vectors fields $X,f_u$, that is 
$$
\Lcal_X f_u(\ubar, \cdot) := \big[ f_u(\ubar, \cdot), X  \big], \qquad \ubar \in \RR. 
$$  
\end{proposition} 

The first term in the right-hand side of \eqref{BV.1} 
depends on the gradient of the unknown function $u$ while the second term depends 
on the function itself. 
Observe that our formula is non-trivial even when $\MM$ is the 
Euclidian space $\RR^d$. In that case, the equation \eqref{BV.1} becomes 
$$
\del_t w + \nabla  (w f_u(u,x)) = - \nabla u \cdot [X, f_u](u,x) - X \cdot \nabla ( \nabla f)(u,x),
$$
where $w = X\cdot \nabla u$. When the flux is independent of $x$ it is obvious that 
for any \emph{constant} vector field the right-hand side of \eqref{BV.1} vanishes and we can recover the 
total variation diminishing property for scalar conservation laws in $\RR^d$. 

We denote by $TV$ the total variation functional for functions defined on the Riemannian manifold $(M,g)$. 
We introduce here a generalization, which extends to discontinuous functions the formula, valid for smooth $u$,  
$$
TV_X(u) := \int_M |X(u(t,x))| \, dv_g. 
$$

\begin{definition} 
\label{BV-2} 

1. To any function $u$ and vector field $X$ defined on the manifold $M$ one associates the quantity 
$$
TV_X(u) := \sup_\phi \int_M u \, \nablac ( \phi \, X) \, dv_g,  
$$
where $dv_g$ is the volume element on the manifold and 
the supremum is taken over all smooth functions $\phi: M \to \RR$ satisfying $\|\phi\|_{L^\infty} \leq 1$. 
When $TV_X(u) < \infty$ the function $u$ is said to have {\bf bounded total variation along $X$.} 

2. A flux $f=f(\ubar, \cdot)$ on $M$ is said to be {\bf divergence-free} if
$$
\nablac f(\ubar, \cdot) = 0, \qquad \ubar\in\RR. 
$$
\end{definition}

\

It follows from Proposition~\ref{BV-1} and the existence theory in \cite{BenArtziLeFloch} that~: 

\begin{theorem} {\rm (Total variation estimates.)} 
\label{BV-3}
Let $u$ be an entropy solution of the conservation law \eqref{IN.1} on the manifold $M$ where $f$ is a smooth flux.  

1. The total variation of $u$ along a vector field $X$ satisfies
\be
\aligned
TV_X(u(t)) \leq TV_X(u(0)) 
& + \sup_{(0,t) \times M} |(\Lie_X f_u)(u, \cdot)|_g \, \int_0^t TV(u(, \cdot)) \, dt
\\
& + \|X(\nablac f)(u,\cdot)\|_{L^1((0,t) \times M)}. 
\endaligned 
\label{Variation1}
\ee 

2. Hence, when $f$ is divergence-free and the Lie bracket of $f_u$ and $X$ vanishes,
$$
\Lcal_X f_u(\ubar, \cdot) = 0, \qquad \ubar \in \RR, 
$$ 
then the solution has bounded total variation along $X$ for all times $t>0$, if this property holds
at the time $t=0$, and moreover,
\be 
TV_X(u(t)) \leq TV_X(u(0)), \qquad t \geq 0. 
\label{BV.2}
\ee  
\end{theorem}

At this juncture, it is important to recall from \cite{BenArtziLeFloch}  that entropy solutions of \eqref{IN.1} 
have bounded variation (in the standard sense) at every time $t$, if this is true at the time $t=0$. 
Actually, one has the general estimate 
$$
TV(u(t)) \le e^{C_1 \, t} \, TV(u(0)) + C_2, 
$$
where $C_1,C_2$ depend on the metric $g$ and on derivatives (up to second order) of the flux $f$. 
In contrast, our result in Theorem~\ref{BV-3} solely assumes that $TV_X(u(0))$ is finite 
and concludes that $TV_X(u(t))$ is finite for all times. 
The property \eqref{BV.2} should be useful when designing approximation schemes for \eqref{IN.1}, 
and it is natural to require \eqref{BV.2} to hold at every time-step of the discretization.

Let us now recall some basic notions and notations that will be used throughout this paper. 
The differential of a function $u : \MM \to \RR$ is the field of $1$-forms defined by 
$$
du(X ) := X (u),  
$$
for every vector field $X$. It is convenient to choose a basis of 
the tangent space $T_x M$ at a point $x \in M$ together with its dual basis associated 
with the cotangent space $T^\star_x M$, so that we can introduce coordinates for 
vectors and covectors relative to these bases. As usual, to each vector $(X^i) \in T_xM$ we associate 
its covector $(g_{ij} \, X^i)$ by lowering indices using the metric (implicit summation on repeated indices being used). 
Denoting $(g^{ij})$ the components of the inverse matrix associated with $(g_{ij})$, 
then, to each covector $(\eta_j) \in \coTM$ we can also associate the vector $(g^{ij} \, \eta_j)$. 

Recall also that the gradient of a (smooth) function $u: M \to \RR$ is the vector field $\nabla_g u$
associated with the field of differential forms $du$. The divergence of a (smooth) vector field $f$ on $M$ 
is the function $\nablac f : M \to \RR$ defined by 
$$
\int_M u \, \nablac f \, dv := -\int_{M}  du (f) \, dv_g 
$$ 
for all smooth $u: M \to \RR$, where $dv_g$ is the volume element on $M$.  

Consider now a coordinate chart $(x^i)$ together with the associated basis of vectors $(\del/\del x^i)$
and covectors $(dx^i)$. The differential of a function $u$ is given by
$$
du  = {\del u \over \del x^i} \, d x^i.
$$ 
In view of the relations $\nabla_g^i u := (\nabla_g u)^i = g^{ij} \, (du)_j$, 
we see that the components of $\nabla_g u$ in the basis $\del/\del x^i$ of $\TM$ are 
$$
(\nabla_g u)^i = g^{ij} \, {\del u \over \del x^j}.
$$
If the support of a function $u$ is included in the domain of definition of the coordinate chart $(x^i)$, 
recalling that $|g|$ denotes the determinant of the metric tensor we can write 
$$ 
\begin{aligned}
\int_{M}  du(f) \, dv_g 
& =\int_{\RR^d} \Big( {\del u \over \del x^i} \, f^i \Big) \, \sqrt{|g|} \, dx^1\cdots dx^d
   \\
& = - \int_{\RR^d} {u \over \sqrt{|g(x)|}} \, 
	 {\del \over \del x^i} \Bigl( f^i \, \sqrt{|g|}\Bigr) \, \sqrt{|g|} \, dx^1\cdots dx^d
\\
& = -\int_M \frac{u}{\sqrt{|g|}}	 {\del \over \del x^i} \Bigl( f^i \, \sqrt{|g|}\Bigr)  \, dv_g, 
\end{aligned}
$$
so that 
\be
\nablac f := {1 \over \sqrt{|g|}} \, 
 {\del \over \del x^i} \Bigl( f^i \, \sqrt{|g|}\Bigr) . 
\label{BV.0.1}
\ee 
 
The expression of $\nablac f$ follows also from computing the
covariant derivatives of the vector field, as follows. 
Given a local chart $(x^i)$ for $M$, the covariant derivative of the vector field $f$ is the $(1,1)$-tensor field 
$\nabla_k f^j$ defined by
$$
\nabla_kf^j := \del_kf^j+\Gamma_{kl}^jf^l,
$$
where $\Gamma_{kl}^j$ are the \emph{Christoffel symbols} given by 
$$
\Gamma_{kl}^j := \frac{1}{2}\big(\del_kg_{li}+\del_lg_{ki}-\del_ig_{kl}\big)g^{ij}.
$$
In particular, we have 
$$
\Gamma_{kj}^j=\frac{1}{2}g^{ij}\del_kg_{ij}
$$
and, since $\del_i|g|=|g|g^{kl}\del_ig_{kl}$, 
$$
\Gamma_{ik}^k=\frac{1}{2|g|}\del_i|g|=\del_i\log\sqrt{|g|}=\frac{1}{\sqrt{|g|}}\del_i\big(\sqrt{|g|}\big).
$$
The divergence of a vector field $f$ is the trace of the covariant derivative, i.e. 
$$
\begin{aligned}
\nablac f := \nabla_kf^k
=\del_kf^k+\Gamma_{ik}^kf^i
& = \del_kf^k+\frac{1}{\sqrt{|g|}}\del_i\big(\sqrt{|g|}\big)f^i
\\
& = \del_kf^k+\frac{1}{\sqrt{|g|}}\big(-\sqrt{|g|}\del_kf^k+\del_k\big(\sqrt{|g|}f^k\big)\big)
\\
& = \frac{1}{\sqrt{|g|}} \, \del_k\big(\sqrt{|g|}f^k\big),
\end{aligned}
$$ 
which is \eqref{BV.0.1}. 


\subsection{First approach} 

The proof of Proposition~\ref{BV-1} will follow from the following two technical lemmas.

\begin{lemma}
\label{BV-4} 
For every smooth function $u: M \to \RR$ the following identities holds ($x \in M$)~: 
\be
\nablac (f(u,x)) = du (f_u (u,x)) + (\nablac f) (u,x)
\label{BV.3}
\ee  
\be
\Lie_X (Y(u,x)) = X(u) Y_u(u,x) + (\Lie_X Y) (u,x)  \qquad \text{ for all vector fields } Y=Y(u,x).
\label{BV.4}
\ee 
\be
X(h(u,x)) = X(u) h_u (u,x) + X(h) (u,x) \qquad \text{ for all functions } h.
\label{BV.5}
\ee 
\end{lemma}

\begin{proof} In local coordinates we have
$$
\begin{aligned}
\nablac (f(u,x)) & = \frac{1}{\sqrt{|g|}} \del_i \big( \sqrt{|g|} f^i(u,x) \big) 
 \\
& = \frac{1}{\sqrt{|g|}} \del_i \big( \sqrt{|g|} f^i \big)(u,x) 
+
 f_u^i (u,x) \del_i u
\\
& = ( \nablac f )(u,x) + du (f_u(u,x)),
\end{aligned}
$$
which proves \eqref{BV.3}. To prove \eqref{BV.4} we note that
$$
\begin{aligned}
\big( \Lie_X (Y(u,x)) \big)^i & = [X, Y(u,x)]^i 
=
 X^j \del_j (Y^i(u,x)) - Y^j(u,x) \del_j X^i
 \\
& =\big( X^j \del_j Y^i - Y^j \del_j X^i \big)(u,x) 
+
 X^j Y^i_u(u,x) \del_j u 
 \\
& = [X, Y]^i (u,x) + X^j \del_j u Y^i_u(u,x). 
\\
\end{aligned}
$$
This shows \eqref{BV.4}. The proof of \eqref{BV.5} is completely similar. 
\end{proof}

\begin{lemma}
\label{BV-5} 
For any (smooth) function $u: M \to \RR$, vector field $X$, and 
flux $f=f (\ubar,\cdot)$, the following identity holds~: 
\be
\begin{aligned}
X \big( \nablac (f(u,x)) \big) 
& =\nablac \big(X(u) f_u(u,x) \big) +  g \big(\nabla_g u,(\Lie_X f_u)(u,x)\big) + X (\nablac f)(u,x). 
\end{aligned}
\label{BV.6}
\ee
\end{lemma}
 
\begin{proof} 
Using \eqref{BV.3} and \eqref{BV.5} we find 
$$
\begin{aligned}
X \big( \nabla_g (f(u,x)) \big) & = X \big( (\nablac f)(u,x) \big)
+
X \big( g( \nabla_g u, (f_u(u,x)) \big)
\\
& = X( \nablac f)(u,x) + X(u)(\nablac f_u)(u,x) 
+
X \big( g( \nabla_g u, (f_u(u,x)) \big).
 \\
\end{aligned}
$$
On the other hand, 
$$
\begin{aligned}
\nablac \big( X(u) f_u(u,x) \big) & = X(u) \nablac (f_u(u,x)) +
	g(\nabla_g (X(u)), f_u(u,x))
\\
	& = X(u) (\nablac f_u)(u,x)+
	X(u) g(\nabla_g u, f_{uu}(u,x)) + g( \nabla_g (X(u)), f_u(u,x)),
\end{aligned}
$$
using \eqref{BV.3}.
Thus, we arrive at 
$$ 
\begin{aligned}
X \big( \nablac (f(u,x)) \big) - \nablac \big( X(u) f_u(u,x) \big) &
 = X(\nablac f)(u,x) + X \big( g( \nabla_g u, f_u(u,x) \big)
\\
& \qquad - X(u) \, g(\nabla_g u, f_{uu}(u,x)) - g(\nabla_g (X(u)), f_u(u,x)).
\end{aligned}
$$ 
Expressing the Leibnitz rule for the Lie derivative we obtain 
$$
X \big( g(\nabla_g u, f_u(u,x)) \big) = (\Lie_X g) \big(\nabla_g u, f_u(u,x) \big) 
+
 g \big(\Lie_X \nabla_g u, f_u(u,x) \big) + g \big(\nabla_g u,\Lie_X (f_u(u,x)) \big).
$$
Therefore, we have 
$$
\begin{aligned}
X \big( \nablac (f(u,x)) \big) - \nablac \big( X(u) f_u(u,x) \big) & 
= X(\nablac f)(u,x) + (\Lie_X g) \big(\nabla_g u, f_u(u,x)\big)
\\
&\quad + g\big( \Lie_X \nabla_g u - \nabla_g (X(u)), f_u(u,x)\big)
\\
&\quad + g \big( \nabla_g u, \Lie_X (f_u(u,x)) - X(u) f_{uu}(u,x))\big).
\end{aligned}
$$ 
Using \eqref{BV.4}, the last term above equals $g \big( \nabla_g u, (\Lie_X f_u)(u,x)\big)$.

The proof of the lemma will be complete once we check that, for all vector fields $X, Z$ and functions $u$,
\be
g \big(\Lie_X \nabla_g u - \nabla_g(X(u)), Z \big) = - (\Lie_X g)(\nabla_g u, Z).
\label{BV.7}
\ee
To this end, in local coordinates we can write 
$$
\begin{aligned}
\big(\Lie_X \nabla_g u - \nabla_g (X(u))\big)^i & = 
X^j \, \del_j \nabla^i u - \nabla^j u \, \del_j X^i - \nabla^i (X^j \, \del_j u)
\\
& = X^j \, \big( \del_j \nabla^i u - \nabla^i \del_j u \bigr) 
                 - \nabla^j u \, \del_j X^i - \nabla^i X^j \, \del_j u. 
\end{aligned}
$$
For the first term above we have 
$$
\begin{aligned}
 X^j \, \big( \del_j \nabla^i u - \nabla^i \del_j u \bigr) &=
 	X^j\, \big(\del_j(g^{ik}) \del_k u + g^{ik}(\del_j \del_k u - \del_k \del_j u)\big)
\\
	&= X(g^{ik})g_{kl} \nabla^l u 
\\
	&= -X(g_{kl})g^{ik}\nabla^l u.
\end{aligned}
$$

For the remaining two terms, a tedious but straightforward computation gives
$$
 - \nabla^j u \, \del_j X^i - \nabla^i X^j \, \del_j u = 
	 - \nabla^j u\, 	g^{il} \big( (\Lie_X g)_{jl} - X(g_{kj}) \big),
$$
where the components of $\Lie_X g$ are computed from the formula
$$
X\big( g(Y, Z) \big) = (\Lie_X g)(Y, Z) + g(\Lie_X Y, Z) + g(Y, \Lie_X Z).
$$
Therefore, we find 
$$
\big(\Lie_X \nabla_g u - \nabla_g (X(u))\big)^i =  - g^{ki} \nabla^j u (\Lie_X g)_{jk}, 
$$
and, finally,
$$
\begin{aligned}
g \big(\Lie_X \nabla_g u - \nabla_g (X(u)), Z \big) &=
 - g_{im} g^{ki} \nabla^j u (\Lie_X g)_{jk} Z^m
 \\
 &= - \nabla^j u (\Lie_X g)_{jk} Z^k
 \\
 &= - (\Lie_X g)(\nabla_g u, Z),
\end{aligned}
$$
which proves \eqref{BV.7}.
This completes the proof of Lemma~\ref{BV-5}. 
\end{proof}

\begin{proof}[Proof of Proposition~\ref{BV-1}]
Applying the vector field $X$ to the conservation law \eqref{IN.1}, we find
$$
\del_t X(u) + X \big( \nablac (f(u, x)) \big) = 0.
$$
The first term is precisely $\del_t w$. Using Lemma~\ref{BV-5} to evaluate the second term, we find 
the desired identity \eqref{BV.1}. 
This completes the proof.
\end{proof}


\subsection{Second approach} 

We provide here a second approach to Proposition~\ref{BV-1}.  
To any vector field $X$ we associate the (local) one-parameter group of diffeomorphisms 
${\varphi : M \to M}$ (defined for all sufficiently small $s$)  from the integral curves of the vector field $X$, 
\be
\frac{d \varphi_s ^i}{ds} (x) = X^i(\varphi_s (x)), \qquad \varphi_0(x) = x. 
\label{BV.8}
\ee 
Given a function $u: M \to \RR$ we may consider the composite function $u \circ \varphi_s$.  

\begin{lemma}
\label{BV-6}
Let $X$ be a smooth vector field on $M$ and 
$(\varphi_s)$ its associated one-parameter group of diffeomorphisms. 
Then, for any smooth function $u : M \to \RR$ it holds 
\be
\aligned 
& \nablac \big( f(u \circ \varphi_s, \cdot)  \big) 
\\
& = \nablac \big( f(u,\cdot)\big) \circ \varphi_s
    - s \, \Big( X(\nablac f) (u, \cdot) + g \big(\nabla_g u, (\Lie_X f_u)(u,\cdot) \big) \Big) + \Ocal(s^2).
\endaligned 
\label{BV.9}
\ee 
\end{lemma}

\begin{proof} We will first check that for any vector field $Y$, 
\be
\label{BV.11}
Y^j(\varphi_s (x)) - \del_i \varphi_s^j Y^i (x) = s \, \Lie_X Y^j + \Ocal(s^2).
\ee
First, it is obvious that 
$$
Y^j \circ \varphi_s = Y^j + s \, X^k \del_k Y^j + \Ocal(s^2).
$$
Second, observe that for all $i,j$ 
$$
\frac{d}{ds} \del_i\varphi_s^j = \del_i \frac{d}{ds} \varphi_s^j = 
\del_i(X^j \circ \varphi_s) = (\del_k X^j) \circ \varphi_s \, \del_i \varphi_s^k.
$$
Also, since $\varphi_0 (x) = x$, we have $\del_i \varphi_s^j\restrict{s=0} = \delta_i^j$
 and so $\frac{d}{ds} \del_i\varphi_s^j\restrict{s=0} = \del_i X^j$.
Therefore, Taylor expanding yields us 
$$
\del_i \varphi^j_s = \delta_i^j + s\, \del_i X^j + \Ocal(s^2).
$$
In consequence we find
$$
\aligned 
Y^j(\varphi_s (x)) - \del_i \varphi_s^j Y^i (x) 
& = s \, X^i \del_i Y^j - s \, Y^i\del_i X^j + \Ocal(s^2) 
\\
& = s \, \Lie_X Y^j + \Ocal(s^2),
\endaligned 
$$
which is precisely \eqref{BV.11}.

Setting $v:=u\circ \varphi_s$, we start from \eqref{BV.3} which, for the function $v$, reads  
\be
\label{BV.10}
\nablac \big( f(v, \cdot)\big) = (\nablac f)(v,\cdot) + dv (f_u(v,\cdot)).  
\ee
Let us compute the last term of the above identity. In local coordinates, we can write 
$$
\aligned 
(dv)_i = \del_i v & = \del_i (u\circ \varphi_s)
\\
& = \big( (\del_j u) \circ \varphi_s \big) \, \del_i  \varphi_s^j, 
\endaligned 
$$
and therefore 
$$
\begin{aligned}
dv (f_u(v,\cdot)) &= \big( (\del_j u) \circ \varphi_s \big) \, \del_i  \varphi_s^j \, f^i_u (v, \cdot)
\\
& = \big( (\del_j u) \circ \varphi_s \big) \, \Big( \del_i  \varphi_s^j \, f^i_u (v, \cdot)
	- f^j_u(u, \cdot) \circ \varphi_s \Big) + \big( \del_j u \, f^j_u(u,\cdot) \big) \circ \varphi_s
\\
& = \big( (\del_j u) \circ \varphi_s \big) \, \Big( \del_i  \varphi_s^j \, f^i_u (v, \cdot)
	- f^j_u(u,\cdot) \circ \varphi_s \Big) +du( f_u(u,\cdot)) \circ \varphi_s.
\end{aligned}
$$

Taking \eqref{BV.11} into account and noting that 
$$
\aligned 
\big( (\del_j u) \circ \varphi_s \big) \, \Lie_X f_u^j(v, \cdot) 
& = \del_j u \,\Lie_X f_u^j(u, \cdot) + \Ocal(s)
\\
&  = g(\nabla_g u, \Lie_X f_u(u,\cdot)) + \Ocal(s), 
\endaligned 
$$
we arrive at  
$$
\begin{aligned}
& dv (f_u(v,\cdot)) 
\\
& = -s \, (\del_j u) \circ \varphi_s (\Lie_X f_u^j)(v,\cdot) + du(f_u(u, \cdot))\circ \varphi_s + \Ocal(s^2).
\\
& = -s \, g \big(\nabla_g u, (\Lie_X f_u)(u,\cdot) \big) + du(f_u(u, \cdot))\circ \varphi_s + \Ocal(s^2)
\\
&= -s \, g \big(\nabla_g u, (\Lie_X f_u)(u,\cdot) \big) + \nablac (f(u,\cdot)) \circ \varphi_s - (\nablac f)(u,\cdot)\circ \varphi_s + \Ocal(s^2),
\end{aligned}
$$
where we have also used \eqref{BV.3} with the function $u$.

Returning now to \eqref{BV.10}, using the identity for $dv (f_u(v,\cdot))$, 
and re-ordering the terms, we find 
$$
\begin{aligned}
\nablac \big(f (v,\cdot) \big) 
= \nablac (f(u,\cdot))\circ \varphi_s  
& + (\nablac f)(v,\cdot)  - (\nablac f)(u,\cdot) \circ \varphi_s 
\\
& - s\, g \big(\nabla_g u, (\Lie_X f_u)(u,\cdot) \big) 
	+ \Ocal(s^2).
\end{aligned}
$$
Finally, observing that 
$$
\begin{aligned}
 (\nablac f)(v,x) - (\nablac f)(v, \varphi_s(x)) &= 
 	- s\frac{d}{ds} (\nablac f)(v, \varphi_s(x))\restrict{s=0} + \Ocal(s^2)
\\
& = - s\,  X^k \, \del_k (\nablac f)(u, \cdot) + \Ocal(s^2)
\\
& = - s\, X (\nablac f)(u,\cdot) + \Ocal(s^2)
\end{aligned}
$$
concludes the derivation of \eqref{BV.9}. 
\end{proof}

\begin{proof}[Proof of Proposition~\ref{BV-1}] We provide here a second proof of this theorem. 
If $u$ is a smooth solution to  \eqref{IN.1} then, 
observing that $\del_t v = (\del_t u)\circ \varphi_s$, we can rewrite \eqref{IN.1} in the form 
$$
\del_t v + \nablac \big( f(u, x) \big) \circ \varphi_s = 0.
$$
So, using \eqref{BV.9} we get
\be
\del_t v + \nablac \big( f(v, \cdot) \big) 
= -  s \, \Big( X(\nablac f) (v, \cdot) + g \big(\nabla_g u, (\Lie_X f_u)(v,\cdot) \big) \Big) + \Ocal(s^2).
\label{BV.11.1}
\ee
Now, recalling the conservation law \eqref{IN.1} and the equation \eqref{BV.11.1}, 
we obtain
\be
\label{BV.12}
\aligned 
& \del_t \Big( \frac{1}{s}(v - u)\Big) + \nablac \Big( \frac{1}{s}(f(v,x) - f(u, x)) \Big)
\\
& = -  X(\nablac f)(v, x) - g(\nabla_g u, (\Lie_X f_u)(v,x)) + \Ocal(s), 
\endaligned
\ee 
in which 
$$
\lim_{s \to 0} \frac{1}{s}\, (v-u)  = X(u) 
$$
and
$$
\aligned
\lim_{s \to 0}  \frac{1}{s}(f(v,\cdot) - f(u, \cdot)) 
&  
= \lim_{s \to 0} \frac{d}{ds} \big( f(v, \cdot) \big) 
\\
& =  f_u(u,\cdot)\, X^i \, \del_i u = X(u) \, f_u(u,\cdot).
\endaligned
$$
So, letting $s$ tend to zero in \eqref{BV.12} yields the desired identity \eqref{BV.1}. 
\end{proof}


\section{Conservation laws on the sphere $S^2$}   
\label{ES-0}

Following \cite{BenArtziLeFloch} 
we include here a discussion of basic properties of the conservation law \eqref{IN.1}
in the case (so important in the geophysical applications) that the manifold $M$ is the sphere $S^2$.  
We discuss here the choice and properties of the flux in \eqref{IN.1} and we emphasize differences between 
the intrinsic and the embedded representations of the sphere. Certain properties will be closely related 
to the fact that the sphere $S^2$ can be embedded isometrically in the Euclidian space $\RR^3$.

\subsection{Spherical coordinates on $S^2$}
One convenient approach consists in parameterizing  the smooth manifold $\Sbold^2$ punctured at the 
North and South poles, 
by using a single chart in spherical coordinates 
$(\varphi,\theta)$, where $\varphi \in [0, 2 \pi]$ denotes the longitude and 
$\theta \in (0, \pi)$ the latitude. Set 
$$
\Omega := T^1\times (0, \pi),
$$ 
where $T^1 = [0,2\pi]$ with the identification of the points $\varphi=0$ and $2 \pi$.  
By definition, a smooth function defined on $\Omega$ is extendable to a smooth
periodic function for all $\varphi \in \RR$. Indeed, in all what follows, we always 
consider smooth functions defined on $\Sbold^2$ represented as functions of
$(\varphi, \theta)$ outside the poles and extended to the poles by continuity. 
We imbed $\Sbold^2$ in  $\RR^3$ according to
$$
x=x(\varphi,\theta)  = \bigl( \cos \varphi \,  \sin \theta, \sin \varphi \,  \sin \theta, \cos \theta \bigr).
$$
The Riemannian metric tensor of the punctured sphere, when expressed in spherical coordinates, read  
$$
\bigl(g_{ij} \bigr) = \begin{pmatrix}
\sin^2 \theta & 0 
\\
0 & 1
\end{pmatrix}. 
$$  
Recall that the inverse matrix of $g=(g_{ij})$ is denoted by $(g^{ij})$. 

The components of the gradient $\nabla_g u$ of a function $u :\Sbold^2 \to \RR$ 
in the basis $(\del / \del\varphi, \del /\del\theta)$ of $\TS$ are
$$
(\nabla_g u)^\varphi = {1\over\sin^2\theta} \, {\del u \over \del \varphi}, 
\qquad 
(\nabla_g u)^\theta = {\del u \over \del \theta}.
$$
According to \eqref{BV.0.1}, the divergence operator on the sphere is given by 
\be
\nabla_g \cdot f := {1 \over \sin \theta} \, 
\Bigl( {\del \over \del \varphi} \bigl( f^\varphi \, \sin \theta\bigr)  
+ 
{\del  \over \del \theta}\bigl( f^\theta \, \sin \theta\bigr)   \Bigr),
\label{ES.1}
\ee
for every smooth vector field $f$ on $S^2$, where $f^\varphi$ and $f^\theta$
are the coordinates of the vector $f$ in the basis $\bigl(\del/\del\varphi, \del/\del\theta\bigr)$. 
Similarly, the Laplace-Beltrami operator is defined 
for all $u: \Sbold^2 \to \RR$ by 
\be 
\Delta_g u := \nabla_g \cdot (\nabla_g u) 
= {1 \over \sin \theta} \, \biggl(
{\del \over \del \varphi} \Bigl( {\del u \over \sin \theta \, \del \varphi} \Bigr)
+
{\del \over \del \theta} \Bigl( {\del u \over \del \theta} \, \sin \theta\Bigr)  
\biggr). 
\label{ES.2}
\ee


We now turn to the following conservation law on the sphere~:  
\be
\del_t u + \nabla_g \cdot f(u, \cdot) = 0, \quad (t,x) \in \RR_+ \times \Sbold^2, 
\label{ES.3}
\ee 
where the flux $f: \RR\times\Sbold^2 \to T\Sbold^2$ is a smooth vector field. 
In view of \eqref{ES.1}, when spherical coordinates are used the conservation law \eqref{ES.3} takes the 
form
\be
\del_t(u\,\sin\theta) + {\del \over \del \varphi} (f^\varphi(u,\varphi,\theta) \, \sin\theta) 
+ {\del \over \del \theta} (f^\theta (u,\varphi,\theta) \, \sin\theta) = 0.
\label{ES.4}
\ee
A periodic boundary condition is imposed in the variable $\varphi \in [0,2 \pi]$, that is, 
\be
u(t, 0,\theta) = u(t,2 \pi, \theta), \quad \theta \in [0, \pi], \, t \geq 0. 
\label{ES.5} 
\ee
On the other hand, no boundary condition may be imposed along $\theta=0$ and  
$\theta = \pi$. Observe that the equation \eqref{ES.4} is \emph{singular}
($\sin \theta = 0$) precisely at the poles.

If $u$ is a smooth solution of the equation \eqref{ES.4}, a natural question is whether $u$ satisfies some additional conservation laws. 
We have already noticed that the Total Variation Diminishing property is available
for solutions to conservation laws in $\RR^d$ when the flux-function is independent of the spatial variable. 
We have seen also that, when writing the conservation law in local coordinates, the flux in general depends
explicitly on the space variables. It would be meaningless to require that the flux is ``independent of $x$''. 
A natural approach is to restrict 
attention to flux that are divergence-free, which in spherical coordinates is expressed as
\be
{\del \over \del \varphi} \big( f^\varphi(\ubar,\varphi,\theta) \, \sin\theta \big) 
+ 
{\del \over \del \theta} \big( f^\theta(\ubar,\varphi,\theta) \, \sin\theta \big) = 0.
\label{ES.6}
\ee
Recall that an entropy / entropy flux pair for the conservation law \eqref{ES.1} is a 
pair $(U,F)$ in which the function $U:\RR\to\RR$ is arbitrary and 
$F(u, x)=\int^u U'(\ubar) \, {\del \over \del u} f(\ubar, x) \, d\ubar$. 

For the sake of completeness we recall~: 

\begin{proposition} {\rm (Intrinsic description.)} 
\label{ES-2}
Consider a conservation law \eqref{IN.1} on the sphere $S^2$ with a divergence-free vector field $f=f(u,x)$, and
let $(U,F)=(U(u),F(u,x))$ be any entropy / entropy flux pair.   
Then, smooth solutions $u=u(t,x)$ of \eqref{ES.4}-\eqref{ES.5} satisfy the additional conservation law  
\be
\del_t U(u) + \nablac F(u, \cdot) = 0. 
\label{ES.8}
\ee
In particular
\be 
{d \over dt} 
\iint_\Omega U(u(t, \cdot))  \, dv_g = 0, 
\quad t \geq 0. 
\label{ES.9}
\ee
\end{proposition}

Of course, for weak solutions, \eqref{ES.8}-\eqref{ES.9} will hold as \emph{inequalities,} only. 
In spherical coordinates, equation \eqref{ES.8} reads
$$
\del_t(U(u)\,\sin\theta) + {\del \over \del \varphi} 
	( F^\varphi(u,\varphi,\theta) \, \sin\theta ) 
+ {\del \over \del \theta} ( F^\theta(u,\varphi,\theta) \, \sin\theta ) = 0.
$$

\begin{proof}
Following \cite{BenArtziLeFloch} we multiply \eqref{IN.1} by $U'(u)$ and get
$$
\del_t (U(u)) + U'(u) \nablac \big( f(u(t, x), x) \big) = 0.
$$
Using \eqref{BV.3} with $F$ instead of $f$, we find
$$
\begin{aligned}
\nablac \big( F(u(t, x), x) \big) &= (\nablac F) (u,x) + du (F_u(u, x))
\\
&= \int_0^u U'(v) (\nablac f_u)(v, x)\, dv + U'(u) du(f_u(u, x)).
\end{aligned}
$$
Again by \eqref{BV.3}, since $f$ is divergence-free, the first term above vanishes and 
$$
U'(u) du(f_u (u,x)) = U'(u) \nablac (f(u, x)),
$$
so that
$$
\nablac \big( F(u(t, x), x) \big) = U'(u) \nablac (f(u, x)),
$$
which proves \eqref{ES.8}.

Furthermore, since $\Sbold^2$ has no boundary, the integral of the divergence of any vector field
vanishes, thus yielding \eqref{ES.9}.

It is illustrative to also derive the equation \eqref{ES.8} using spherical coordinates.
To simplify the notation we set  
$$
h^i = f^i\sin\theta, \quad i=\varphi, \theta; 
\qquad 
\del_\alpha = \frac{\del}{\del\alpha}, \quad \alpha = \varphi,\theta.
$$  
We denote here by $h^i_u, h^i_\varphi, h^i_\theta$ respectively
 the derivatives of $h^i$ with respect to the first, the second, and the third variables. 
Then, the equation \eqref{ES.8} becomes
$$
\del_t(u\,\sin\theta)+\del _\varphi\big(h^\varphi(u,\varphi,\theta)\big)+
\del_\theta\big(h^\theta(u,\varphi,\theta)\big)=0.
$$
Multiplying it by $U'(u)$, we get
$$
\begin{aligned}
0&=\del_t(U(u)\,\sin\theta)+U'(u)\del_\varphi(h^\varphi)+U'(u)\del_\theta(h^\theta)\\
&=\del_t(U(u)\,\sin\theta)+U'(u)\big(h^\varphi_\varphi+h^\varphi_u\del\varphi u\big)+U'(u)\big(h^\theta_\theta+h^\theta_u\del_\theta u\big),\\
\end{aligned}
$$
or equivalently, with $H:=F\sin\theta$, so that $H^i_u=F^i_u\sin\theta$,
$$
\begin{aligned}
&\del_t(U(u)\,\sin\theta)+H^\varphi_u\del_\varphi u+H^\theta_u\del_\theta u+U'(u)(h^\varphi_\varphi+h^\theta_\theta)\\
&=\del_t(U(u)\,\sin\theta)+\del_\varphi H^\varphi+\del_\theta H^\theta+Q=0,
\end{aligned}
$$
where
$$
Q=U'(u)(h^\varphi_\varphi+h^\theta_\theta)-H^\varphi_\varphi-H^\theta_\theta.
$$
Suppose now that the condition \eqref{ES.6} is satisfied, that is, 
$h^\varphi_\varphi+h^\theta_\theta=0$ for all $u\in\RR$. Then, we have 
$$
\begin{aligned}
H^\varphi_\varphi+H^\theta_\theta&=\int_0^uU'(v)(h^\varphi_{u,\varphi}+h^\theta_{u,\theta}) \, dv
&=\int_0^uU'(v)(h^\varphi_{\varphi}+h^\theta_{\theta})_u\,dv=0,
\end{aligned}
$$
and therefore $Q=0$, which completes the derivation of \eqref{ES.8} in spherical coordinates. 
\end{proof}


\subsection{Embedded representation of $S^2$}
 
In the embedded approach, the sphere $\Sbold^2$ is regarded as a submanifold of the Euclidian space $\RR^3$, 
via some  isometric embedding $h: \Sbold^2 \to \RR^3$.
The canonical Riemannian metric $g$ on $\Sbold^2$ is  induced by the canonical metric $\widehat g$
on the Euclidian space $\RR^3$, by $g = h^\star\widehat g$, where, 
by definition, the \emph{pull-back} $h^\star$ maps the $2$-tensor field $\widehat g$ 
on $\RR^3$ to the $2$-tensor field $g$ on $\Sbold^2$, as follows:  
\be
(h^\star \widehat g)(\xi,\eta) := \widehat g(h_\star\xi, h_\star\eta), 
\qquad 
\xi,\eta\in T\Sbold^2. 
\label{ES.10}
\ee 
The tangent vector $h_\star\xi$ is the \emph{push-forward} of $\xi\in T\Sbold^2$ into $\RR^3$:
\be 
(h_\star\xi)(f) := \xi( f \circ h) \quad \text{ for every function } f : \RR^3 \to \RR. 
\label{ES.11}
\ee 

From now on we fix some local coordinates $(x^i)_{i=1,2}$ on $\Sbold^2$ (for instance, the spherical 
coordinates described above) as well as
some orthonormal coordinates $(\widehat x^\alpha)_{\alpha=1,2,3}$ on $\RR^3$.
In particular, we have $\widehat g_{\alpha\beta} = \delta_{\alpha\beta}$ (the Kronecker symbol). 

We have
\be
\label{ES.12}
(h_\star\xi)^\alpha=\frac{\del \widehat x^\alpha}{\del x^i}\xi^i,
\ee
so that we may compute the expression of the metric tensor $g$ in these local coordinates:
\be 
\begin{split}
g_{ij}\xi^i\eta^j
& = \widehat g_{\alpha\beta} \, (h_\star\xi)^\alpha \, (h_\star\eta)^\beta 
\\
& = {\del \widehat x^\alpha \over \del x^i} \, {\del \widehat x_\alpha \over \del x^j} \, \xi^i \, \eta^j,
\end{split}
\nonumber
\ee 
where $\widehat x_\alpha := \widehat g_{\alpha \beta}\widehat x^\beta$. 
From this, it follows that 
\be
g_{ij} = {\del \widehat x^\alpha \over \del x^i} \, {\del \widehat x_\alpha \over \del x^j}.
\label{ES.13}
\ee 
In other words, we have $(g_{ij}) = J^t \, J$, where $J= \big( \del \widehat x^\alpha / \del x^j\big)$ is 
the Jacobian matrix of the embedding $h$.

In what follows, all vectors are considered to be in $\RR^3$, unless otherwise stated. 
A vector tangent to the submanifold $\Sbold^2$ of $\RR^3$ is regarded as an element of $\RR^3$ through the identification
$$
\TS \subset T_x\RR^3 \simeq \RR^3.
$$
Consider now the spherical coordinates introduced earlier.
The outward unit vector to $\Sbold^2$ in $\RR^3$ at the point $x=(\varphi,\theta)$ is   
given by 
$$ 
h(\varphi, \theta) = \nn(x) := \bigl( \sin \theta \, \cos \varphi, 
              \sin \theta \, \sin \varphi, \cos \theta \bigr).
$$
At any given point $x=(\varphi,\theta)$, 
a normalized basis of the tangent plane to $\Sbold^2$ in $\RR^3$ (now viewed as a two-dimensional subspace of $\RR^3$) 
is given by the two vectors 
$$ 
\aligned 
&  \ttt_1(x) :={1\over\sin \theta } {\del \nn \over  \del \varphi} 
      = \bigl( - \sin \varphi, \cos \varphi, 0 \bigr), 
      \\
& \ttt_2(x) := - {\del \nn \over \del \theta}
      =\bigl( - \cos \theta \, \cos \varphi, 
            - \cos \theta \, \sin \varphi, \sin \theta \bigr). 
\endaligned 
$$  

Given a vector $\xi$ in $TS^2$ 
with (intrinsic) spherical coordinates $(\xi^\varphi, \xi^\theta)$, we associate to it the vector $\tilde \xi\in \RR^3$ defined by \eqref{ES.11}.
The components of $\tilde \xi$ in the basis $(\ttt_1,\ttt_2, \nn)$ of $\RR^3$ are denoted 
by $(\tilde \xi^i)_{1\le i \le 3}$. The vector $\tilde\xi$ is tangent to the sphere: using \eqref{ES.12} 
to compute its coordinates, we find
$$
\tilde \xi = \frac{\del \nn}{\del \varphi} \xi^\varphi + 
	\frac{\del \nn}{\del \theta} \xi^\theta = \sin \theta\,\xi^\varphi \ttt_1 -
	\xi^\theta \ttt_2,
$$
so that in spherical coordinates the intrinsic and embedded coordinates of a tangent vector are related by
\be
\label{ES.14}
\tilde\xi^1 = \sin\theta \,\xi^\varphi,  \qquad \tilde\xi^2 = - \xi^\theta, \qquad \tilde \xi^3 = 0.
\ee

Observe that the normal vector is well-defined for all $\theta$ (and $\varphi$). However, 
the North ($\theta= 0$) and South ($\theta=\pi$) poles are 
singularities for the tangent vectors in spherical coordinates. Namely, 
the particular value of $\varphi$ is irrelevant when $\theta= 0$ or $\pi$; however,  
the tangent vectors to $\Sbold^2$ do depend upon $\varphi$ when $\theta= 0$ or $\pi$, 
so that a basis of the tangent space is not uniquely specified at the poles. 
This is consistent with the well-known fact that any globally smooth, tangent 
vector field on the sphere must vanish at least once.


In (embedded) spherical coordinates, 
the flux $f: \RR \times \Omega \to \RR^3$ has the general form 
$$ 
f(u, \varphi,\theta) =: \tilde f^1(u, \varphi,\theta) \, \ttt_1(\varphi,\theta)  
                      + \tilde f^2(u, \varphi,\theta) \, \ttt_2(\varphi,\theta). 
$$
Note that $f$ is tangent to the sphere. In view of \eqref{ES.1}, \eqref{ES.14}, we have 
$$
\tilde f^1 = \sin\theta \, f^\varphi,  \quad \tilde f^2 = - f^\theta,
$$
and the divergence operator takes the form
$$
\nablac f := {1 \over \sin \theta} \, \Bigl( {\del  \tilde f^1\over \del \varphi}  
-  {\del \over \del \theta} \bigl( \tilde f^2 \, \sin \theta\bigr)  
  \Bigr). 
$$ 
In terms of the unknown function $u=u(t,\varphi,\theta)$,  
\eqref{ES.3} is a conservation law with non-constant coefficients
posed in the rectangular domain $(t,\varphi,\theta) \in \RR_+ \times \Omega$: 
\be 
\del_t \bigl( u \, \sin \theta \bigr)  
+ 
{\del \over \del \varphi} \bigl(\tilde f^1(u, \varphi,\theta) \bigr) 
- 
{\del \over \del \theta} \bigl(\tilde f^2(u, \varphi,\theta) \, \sin \theta\bigr)  
= 0. 
\label{ES.15}
\ee

For the sake of clarity we will always state the final results in the intrinsic notation $f^\varphi, f^\theta$, 
but in the course of the calculation it will be convenient to use the embedded notation $\tilde f^1, \tilde f^2$, 
which has the advantage to be based on an orthonormal basis of $\RR^3$. 

Let us now define the tubular neighborhood
 $T \subset \RR^3$ of $\Sbold^2$ consisting of all spheres with radius
$1/2 < r < 3/2$, so that $\Sbold^2$ corresponds to $r=1$. Then, any
smooth tangent vector field $f(u, x)$, $x\in \Sbold^2$, may be represented 
at each point of the sphere by 
\be
f(u,x) = \nn(x) \times \Phi \big(u, \nn (x)\big),
\label{ES.16}
\ee
where $\Phi(u, x)$, $x \in T$ is a smooth vector field in $T \subset \RR^3$.
We can thus write $\Phi(u, x) = \Phi(u, \nn(x)).$

\begin{lemma} 
\label{ES-4}
With the notation \eqref{ES.16} the equation \eqref{ES.15} takes the form 
\be 
\del_t \bigl( u \, \sin \theta \bigr)  
+ 
{\del \over \del \varphi} \Bigl(  \Phi(u,\nn( \varphi, \theta) )  \cdot {\del \nn \over \del \theta}  \Bigr) 
- 
{\del \over \del \theta} \Bigl(  \Phi(u,\nn( \varphi, \theta) )  \cdot {\del \nn \over \del \varphi} \Bigr) 
= 0,
\label{ES.17}
\ee
\end{lemma}

\begin{proof} Setting 
$$
A :=\begin{pmatrix} 
\sin \theta \, \cos \varphi  &  - \sin \varphi &  - \cos \theta \, \cos \varphi   
\\
\sin \theta \, \sin \varphi  &   \cos \varphi  &  - \cos \theta \, \sin \varphi    
\\
\cos \theta                         &     0                         &  \sin \theta     
\end{pmatrix}
$$
we have 
$$
\nn = A \, \begin{pmatrix} 1 \\ 0 \\ 0 \end{pmatrix},
\quad 
\ttt_1 = A \, \begin{pmatrix} 0 \\ 1 \\ 0 \end{pmatrix},
\quad 
\ttt_2 = A \, \begin{pmatrix} 0 \\ 0 \\ 1 \end{pmatrix},
$$ 
and $ A^{-1} = A^t$. Since $\det(A) = 1$ the basis $( \ttt_1, \ttt_2, \nn)$ is positively oriented and
$$
\nn \times \ttt_1 = \ttt_2, \qquad \nn \times \ttt_2 = - \ttt_1.
$$
Using coordinates, we have 
$$
\tilde \Phi^i(u, \nn(\varphi, \theta)):= \Phi(u, \nn(\varphi, \theta)) \cdot \ttt_i(\varphi, \theta)
$$ 
and
$$
\aligned 
f(u, \varphi, \theta) 
& = \nn(\varphi, \theta) \times \Big( 
\tilde \Phi^1(u, \nn(\varphi, \theta)) \, \ttt_1(\varphi, \theta) 
         +\tilde \Phi^2(u, \nn( \varphi, \theta)) \, \ttt_2(\varphi, \theta) \Big) 
\\
& =  
\tilde \Phi^1(u, \nn( \varphi, \theta)) \, \ttt_2( \varphi, \theta) 
- \tilde \Phi^2(u, \nn( \varphi, \theta)) \, \ttt_1(\varphi, \theta),
\endaligned 	 
$$
so that $\tilde f^1 =  - \tilde \Phi^2$, $\tilde f^2 = \tilde \Phi^1.$ But 
$$
\tilde \Phi^1 = \frac{1}{\sin \theta} \Big( \Phi \cdot \frac{\del \nn}{\del \varphi}
 \Big), \quad \tilde \Phi^2 = - \Big( \Phi\cdot \frac{\del \nn}{\del \theta} \Big),
 $$
and so 
$$
\tilde f^1 =\Big( \Phi\cdot \frac{\del \nn}{\del \theta} \Big),
 \quad \tilde f^2 =\frac{1}{\sin \theta} \Big( \Phi \cdot \frac{\del \nn}{\del \varphi}
 \Big), 
$$ 
or equivalently using the intrinsic notation 
$$
 f^\varphi =\frac{1}{\sin \theta} \Big( \Phi\cdot \frac{\del \nn}{\del \theta} \Big),
 \quad  f^\theta = - \frac{1}{\sin \theta} \Big( \Phi \cdot \frac{\del \nn}{\del \varphi}
 \Big), 
$$
which, in view of \eqref{ES.15}, leads to \eqref{ES.17}.
\end{proof}

Observe that, for smooth solutions, \eqref{ES.17} takes the equivalent form 
\be 
\del_t \bigl( u \, \sin \theta \bigr)  
+ \del_u \Phi( u, \nn(\varphi,\theta))  \cdot \Bigl( 
{\del u\over \del \varphi} \, {\del \nn \over \del \theta} 
- {\del u \over \del \theta} \, {\del \nn \over \del \varphi} \Bigr)   
= 0.
\label{ES.18} 
\ee
We have the analogue of Proposition~\ref{ES-2}:

\begin{proposition} {\rm (Embedded description.)} 
\label{ES-5}
Suppose that the vector field $\Phi$ satisfies the divergence-free condition: 
for all $\ubar\in\RR$ and  $(\varphi,\theta) \in \Omega$, 
\be
{\del \over \del \varphi} \Bigl(  \Phi(\ubar, \nn(\varphi,\theta))  \cdot {\del \nn \over \del \theta}(\varphi,\theta)  \Bigr) 
-
{\del \over \del \theta} \Bigl(  \Phi(\ubar, \nn(\varphi,\theta))  \cdot {\del \nn \over \del \varphi}(\varphi,\theta) \Bigr) 
= 0 
\label{ES.19}
\ee
or, equivalently,
\be
{\del \Phi \over \del \varphi} (\ubar) \cdot {\del \nn \over \del \theta}
-
{\del \Phi \over \del \theta} (\ubar) \cdot {\del\nn \over \del \varphi}
= 0. 
\label{ES.20}
\ee 
Let $U:\RR\to\RR$ and $\Psi:\RR\times\Sbold^2\to\RR^3$ be smooth maps such that 
\be
\Psi(u, \nn(\varphi, \theta)) =
 \int_0^u U'(\ubar) \, {\del \over \del u} \Phi(\ubar, \nn(\varphi, \theta)) \, d\ubar.
\label{ES.21}
\ee 
Then every smooth solution $u$ of \eqref{ES.17} satisfies
\be 
\del_t \bigl( U(u) \, \sin \theta \bigr)  +
{\del \over \del \varphi} \Bigl(  \Psi(u, \nn( \varphi, \theta))  
\cdot {\del \nn \over \del \theta}  \Bigr) 
-
{\del \over \del \theta} \Bigl(  \Psi(u, \nn(\varphi, \theta))  \cdot {\del \nn \over \del \varphi} \Bigr) 
= 0. 
\label{ES.22}
\ee 
\end{proposition}

\begin{proof} Multiplying the equation \eqref{ES.17} by $U'(u)$ we obtain 
$$
\del_t \bigl(U(u) \, \sin \theta \bigr)  
- \del_u \Psi(u, \nn( \varphi, \theta)) \cdot \del_\theta u\,\frac{\del\nn}{\del\varphi} 
+ \del_u \Psi(u, \nn(\varphi, \theta)) \cdot \del_\varphi u\,\frac{\del\nn}{\del\theta} =0, 
$$
thus 
$$
\del_t \bigl(U(u) \, \sin \theta \bigr) 
-\del_\theta \Big( \Psi (u, \nn(\varphi, \theta)) \cdot \frac{\del\nn}{\del\varphi} \Big)
+\del_\varphi \Big( \Psi(u, \nn(\varphi, \theta)) 
\cdot \frac{\del\nn}{\del\theta} \Big) + Q = 0,
$$
where
$$
\begin{aligned} 
Q  & : = \Psi_\theta\cdot\frac{\del\nn}{\del\varphi}-\Psi_\varphi\cdot\frac{\del\nn}{\del\theta}
\\
& = \int^u U'(v) \, \Phi_{u\theta}\,dv\frac{\del\nn}{\del\varphi}-\int^uU'(v)\Phi_{u\varphi}\,dv\frac{\del\nn}{\del\theta}
\\
& = \int^uU'(v)\Big(\Phi_{\theta}\cdot\frac{\del\nn}{\del\varphi}-\Phi_\varphi\cdot\frac{\del\nn}{\del\varphi}\Big)_u\,dv=0,
\end{aligned}
$$
by the condition \eqref{ES.20}.
\end{proof}

\begin{remark}
\label{ES-6}
The constraint \eqref{ES.19} can also be written in the form
$$
\Big( D \Phi \cdot \frac{\del \nn}{\del \varphi} \Big) \cdot 
	\frac{\del \nn}{\del \theta} -
\Big( D \Phi \cdot \frac{\del \nn}{\del \theta} \Big) \cdot
	\frac{\del \nn}{\del \varphi} = 0.
$$
In particular, if the map $\Phi$ coincides with the gradient of some
potential $h$, that is, $\Phi = D h$ for $h: T \to \RR$, then 
the divergence-free condition \eqref{ES.19} is automatically satisfied.
\end{remark} 


\section{Total variation diminishing estimates on the sphere $S^2$}   
\label{BVS-0}

We are now in a position to discuss the total variation of solutions 
to scalar conservation laws on the sphere. We are going to investigate the conditions 
discovered in Theorem~\ref{BV-3}. After comparing the intrinsic and embedded 
description it turned out that the second approach led to 
simpler (but still not so straightforward) and more intuitive formulas. 
In view of Section~\ref{ES-0}, a vector field tangent to the embedded sphere $\Sbold^2$ 
can be expressed as the cross product between the exterior normal
$\nn (x)$ and some other vector field which can always be chosen to be tangent to the sphere. 
So, 
given arbitrary flux $f$ and vector field $X$ let us introduce (smooth) maps
$\Phi: \RR \to T \to \RR^3$ and $\Psi: T \to \RR^3$ (which are not unique) such that
\be
\label{TVS.1} 
\begin{aligned}
& f(u, \varphi, \theta) = \nn(\varphi, \theta) \times \Phi(u,  \nn(\varphi, \theta))
\\
& X(\varphi, \theta) = \nn(\varphi, \theta) \times \Psi( \nn(\varphi, \theta))
\end{aligned}
\ee
(where $T\subset \RR^3$ is the tubular neighborhood of $\Sbold^2$ defined in Section~\ref{ES-0}).
It should be kept in mind that all the vectors under consideration here belong to $\RR^3$.

Having in mind the class of divergence-free flux, we introduce~: 

\begin{definition} The flux $f$ and vector field $X$ are said to satisfy the {\bf gradient condition}  
if there exist functions ${a: \RR \times \RR^3 \to \RR}$ and ${b : \RR^3 \to \RR}$ 
depending solely on $\nn=x/|x|$ and are smooth in a neighborhood of the $2$-sphere
such that
\be
\Phi(u,x) = D_x a(u,x/|x|), \qquad \Psi(x) = D_x b(x/|x|), \qquad x \in \RR^3.
\label{BVK}
\ee
\end{definition}

From the gradient condition we deduce that, for instance,  
$$
\Phi = D_x a = D_\nn a \cdot D_x(x/|x|).
$$
Clearly, $D_x(x/|x|)\cdot \nn$ vanishes identically, since it is the derivative of the map $x/|x|$ in the direction of the vector $\nn$. 
Therefore,
$$
\Phi\cdot \nn = \Psi \cdot \nn = 0, 
$$
that is, the vectors $\Phi$ and $\Psi$ are tangent to $\Sbold^2$. Observe also that as in Remark \ref{ES.6},
the divergence-free condition is automatically satisfied for fields satisfying the gradient condition.

Here is the main result of this section.

\begin{theorem}
\label{TVS-1}
Let $f$ and $X$ be a flux and a vector field defined on the sphere $\Sbold^2$ respectively, 
satisfying the gradient condition \eqref{TVS.1}-\eqref{BVK} for some maps $\Phi,\Psi$.  
Then, the flux $f$ is divergence-free and one has 
\be
\label{TVS.2}
\Lie_X f_u (\ubar, \cdot)= 0 \, \text{ if and only if } \, f_u(\ubar, \cdot) = C(\ubar, \cdot) \, X, 
\qquad \ubar \in \RR,  
\ee
for some scalar-valued flux-function $C=C(\ubar,\cdot)$ satisfying $X(C(\ubar, \cdot)) = 0$.

When \eqref{TVS.2} holds, and when $(\nablac f)(u, \cdot) = 0$, the conservation law \eqref{ES.15}
decouples into a one-parameter family of independent, one-dimensional conservation laws. 
\end{theorem}

Hence, in the case of flux which are ``gradients'', we are able to identify the class of flux
satisfying the condition of Section~\ref{BV-0}. For each $\ubar$, the vectors $f_u(\ubar, \cdot)$ and $X$ are parallel 
and $C(\ubar, \cdot)$ is constant along the integral curves of the vector field $X$. The curved geometry 
of the sphere restrict the class of flux enjoying the TVD property. This should be compared with the 
case of equations in the Euclidian space $\RR^d$ where no such condition is implied on the flux.

Throughout the rest of this section, a value $u\in \RR$ is fixed and so we omit it. 
Recall from Section~\ref{ES-0} that the components  $(\tilde f^1, \tilde f^2)$
of $f$ in the basis $(\ttt_1, \ttt_2)$ of the tangent plane to $\Sbold^2$ at a point $x(\varphi, \theta)$
read~:
$$
\hfill\tilde f^1 = - \Phi(\nn) \cdot \nn_\theta, \quad
\tilde f^2 = \frac{1}{\sin \theta} \Phi(\nn)\cdot \nn_\varphi,
$$
where we have used subscripts to indicate differentiation.
Using the equation \eqref{ES.14} which connects these components 
with the intrinsic components $(f_u^\varphi, f_u^\theta)$ of $f_u$, we find
$$ 
\begin{aligned}
& f^\varphi = \frac{1}{\sin \theta} \Phi(\nn) \cdot \nn_\theta
\\
& f^\theta = - \frac{1}{\sin \theta} \Phi(\nn)\cdot \nn_\varphi.
\end{aligned}
$$
A similar formula holds for $X$~:
$$ 
\begin{aligned}
& X^\varphi = \frac{1}{\sin \theta} \Psi(\nn) \cdot \nn_\theta
\\
& X^\theta = - \frac{1}{\sin \theta} \Psi(\nn)\cdot \nn_\varphi.
\end{aligned} 
$$

Observe now the following elementary identities:
\be
\begin{aligned}
&\nn_{\varphi\varphi} = -\sin^2 \theta \, \nn -\sin\theta \cos\theta \, \nn_\theta,
\\
&\nn_{\varphi \theta} = \frac{\cos \theta}{\sin \theta} \nn_\varphi
\\
&\nn_{\theta \theta} = -\nn.
\end{aligned}
\label{BVK.16}
\ee 

We are now in a position to reformulate condition $\Lie_X f_u = 0$.

\begin{lemma}
\label{BVK-7}
The condition $\Lie_X f_u = 0$ is equivalent to the requirement that 
$\Psi=\Psi (\nn)$ satisfies the following system of partial differential equations with unknown $X$~:  
\be 
\begin{aligned}
	&\begin{aligned}
	(\Phi_u \cdot \nn_\theta) (\Psi_\varphi \cdot \nn_\theta)
		- (\Phi_u \cdot \nn_\varphi) (\Psi_\theta \cdot \nn_\theta)
	& = - (\Psi \cdot \nn_\varphi) \big( (\Phi_{u\theta} \cdot \nn_\theta)
		- (\Phi_u \cdot \nn) \big)
	\\
	&\quad - (\Psi \cdot \nn) (\Phi_u \cdot \nn_\varphi) 
		+ (\Psi \cdot \nn_\theta) (\Phi_{u\varphi} \cdot \nn_\theta),
	\end{aligned}
	\\
	&\begin{aligned}
	(\Phi_u \cdot \nn_\theta) (\Psi_\varphi \cdot \nn_\varphi)
		- (\Phi_u \cdot \nn_\varphi) (\Psi_\theta \cdot \nn_\varphi)
	& = (\Psi \cdot \nn_\theta) \big( (\Phi_{u\varphi} \cdot \nn_\varphi)
		- \sin^2 \theta (\Phi_u \cdot \nn) \big)
	\\
	&\quad + \sin^2 \theta (\Psi \cdot \nn) (\Phi_u \cdot \nn_\theta) 
		- (\Psi \cdot \nn_\varphi) (\Phi_{u\theta} \cdot \nn_\varphi).
	\end{aligned}
\end{aligned}
\label{BVK.18}
\ee
\end{lemma}

\begin{proof} The condition on the first component, $(\Lie_X f_u)^\varphi$, of the vector
$\Lie_X f_u$ reads 
$$
 (X^\varphi \del_\varphi + X^\theta \del_\theta)
	f_u^\varphi - (f_u^\varphi \del_\varphi + f_u^\theta \del_\theta) X^\varphi = 0, 
$$
that is
$$
	\begin{aligned}
	&\frac{1}{\sin \theta}(\Psi \cdot \nn_\theta) 
		\del_\varphi(\Phi_u \cdot \nn_\theta)
		-(\Psi \cdot \nn_\varphi) \del_\theta \Big(\frac{1}{\sin\theta}
		(\Phi_u \cdot \nn_\theta)\Big) 
	\\
	&\qquad - \frac{1}{\sin\theta}(\Phi \cdot \nn_\theta)
		\del_\varphi(\Psi \cdot \nn_\theta) + (\Phi \cdot \nn_\varphi)
		\del_\theta \Big(\frac{1}{\sin \theta}(\Psi \cdot \nn_\theta)\Big) =0. 
	\end{aligned}
$$
This equation takes the form 
$$ 
\begin{aligned}
	&(\Psi \cdot \nn_\theta) (\Phi_{u\varphi} \cdot \nn_\theta)
		-(\Psi \cdot \nn_\varphi) (\Phi_{u\varphi} \cdot \nn_\theta)
		 + \frac{\cos \theta}{\sin \theta} (\Psi \cdot \nn_\varphi)
		(\Phi_u \cdot \nn_\theta)
	\\
	&\qquad - (\Phi_u \cdot \nn_\theta)
		(\Psi_\varphi \cdot \nn_\theta) + (\Phi_u \cdot \nn_\varphi)
		(\Psi_\theta \cdot \nn_\theta)  
		- \frac{\cos \theta}{\sin \theta} (\Phi_u \cdot \nn_\varphi)
		(\Psi \cdot \nn_\theta)
	\\
	&\qquad + (\Psi \cdot \nn_\theta) (\Phi_u \cdot \nn_{\varphi\theta})
		- (\Psi \cdot \nn_\varphi) (\Phi_u \cdot \nn_{\theta\theta})
	\\
	&\qquad - (\Phi_u \cdot \nn_\theta) (\Psi \cdot \nn_{\theta\varphi})
		+ (\Phi_u \cdot \nn_\varphi) (\Psi \cdot \nn_{\theta\theta}) =0.
	\end{aligned} 
$$
In view of \eqref{BVK.16}, all terms above which contain the factor $(\cos \theta / \sin \theta)$ cancel out. 
Using that $\nn_{\theta\theta} = -\nn$ and re-arranging the terms leads to the first equation.

The second equation is obtained similarly, by developing 
$$
(\Lie_X f_u)^\theta = (X^\varphi \del_\varphi + X^\theta \del_\theta)
	f_u^\theta - (f_u^\varphi \del_\varphi + f_u^\theta \del_\theta) X^\theta
$$
and using \eqref{BVK.16}. 
\end{proof}

We now show: 

\begin{lemma}
\label{BVK-8}
A point $\overline x \in \Sbold^2$ being fixed, there exist non-trivial flux $f$ and vector fields $X$ 
satisfying the gradient condition, and possibly defined in a small neighborhood of $\overline x$ only, 
such that $\Lie_X f_u = 0$. These vector fields are characterized explicitly by the condition  
$$
X \times f_u = \widetilde C \, \nn
$$
for some scalar $\widetilde C = \widetilde C(u)$ (independent of $x$).
Among these solutions, the only flux $f$ and vector field $X$ that are globally defined and smooth
on the whole sphere $\Sbold^2$ 
are given by  
$$
f_u(u, \cdot) = C(u, \cdot) \, X,
$$
where $C=C(,\cdot)$ is a (scalar-valued) function such that $X(C(u, \cdot)) =0$.  
\end{lemma}

Vector fields that are defined locally on $\Sbold^2$ may still be useful to
 control the total variation in a domain of dependence, but they do not provide information on 
the large-time behavior of solutions.  

\begin{proof}
First, let us show that for $i = \varphi, \theta$ we have
$a_i = \Phi\cdot \nn_i$ and $b_i = \Psi \cdot \nn_i$ (subscripts indicate 
differentiation). A straightforward calculation shows that
$$
D_x a = D_\nn a\cdot D_x \left(\frac{x}{|x|} \right) =
 D_\nn a \cdot \left( \frac{\mathrm{Id} - \nn \otimes \nn}{|x|} \right).
$$
In particular, restricting the formula to the sphere we find $D_x a = D_\nn a \cdot 
(\mathrm{Id} - \nn \otimes \nn)$, and therefore
$$
a_i = D_\nn a \cdot \nn_i = D_\nn a \cdot (\mathrm{Id} - \nn \otimes \nn) 
+ D_\nn a  \cdot (\nn\otimes \nn) \cdot \nn_i = D_x a \cdot \nn_i = \Phi \cdot \nn_i,
$$
because $(\nn \otimes \nn) \cdot \nn_i = 0$. In consequence, we now simply write $Da$, 
$Db$. Recalling that by assumption we have $Da\cdot \nn = Db \cdot \nn = 0$,
 we can rewrite \eqref{BVK.18} in terms of $b$:
 \be
 \begin{aligned}
& a_{u\theta} (\del_\varphi Db \cdot \nn_\theta) 
- a_{u\varphi} (\del_\theta Db \cdot \nn_\theta)
= - b_\varphi (\del_\theta Da_u \cdot \nn_\theta) 
+ b_\theta (\del_\varphi Da_u \cdot \nn_\theta),
\\
&
 a_{u\theta} (\del_\varphi Db \cdot \nn_\varphi) 
- a_{u\varphi} (\del_\theta Db\cdot \nn_\varphi) 
=  b_\theta (\del_\varphi Da_u \cdot \nn_\varphi) 
- b_\varphi (\del_\theta Da_u \cdot \nn_\varphi).
\end{aligned}
 \label{BVK.19}
\ee
This is now a system in the \emph{scalar} variable $b$. 
Furthermore, using \eqref{BVK.16} we find
$$
\begin{aligned}
&a_{u\theta\theta} = \del_\theta (Da_u \cdot \nn_\theta) 
	= \del_\theta Da_u \cdot \nn_\theta - Da_u \cdot \nn =
	\del_\theta Da_u \cdot \nn_\theta,
\\
&a_{u\varphi\varphi} = \del_\varphi Da_u \cdot \nn_\varphi
	+ Da_u \cdot \nn_{\varphi\varphi}
	=  \del_\varphi Da_u \cdot \nn_\varphi - 
	\sin \theta \cos \theta \, a_{u\theta},
\\
&a_{u\theta\varphi} = \del_\varphi Da_u \cdot \nn_\theta 
+ Da_u \cdot \nn_{\varphi\theta} =  \del_\varphi Da_u \cdot \nn_\theta 
+ \frac{\cos \theta}{\sin \theta} Da_u \cdot \nn_\varphi 
= \del_\varphi Da_u \cdot \nn_\theta 
+ \frac{\cos \theta}{\sin \theta} a_{u\varphi},
\end{aligned}
$$
so that
$$
\begin{aligned}
&\del_\theta Da_u \cdot \nn_\theta = a_{u\theta\theta},
\\
&\del_\varphi Da_u \cdot \nn_\varphi = a_{u\varphi\varphi} +
\sin \theta \cos \theta \, a_{u\theta},
\\
&\del_\varphi Da_u \cdot \nn_\theta = a_{u\theta\varphi}
- \frac{\cos \theta}{\sin \theta} a_{u\varphi}.
\end{aligned}
$$
Similarly, we obtain 
$$ 
\begin{aligned}
&\del_\theta Db \cdot \nn_\theta = b_{\theta\theta}
\\
&\del_\varphi Db \cdot \nn_\varphi = b_{\varphi\varphi} +
\sin \theta \cos \theta \, b_\theta
\\
&\del_\varphi Db \cdot \nn_\theta = b_{\theta\varphi}
- \frac{\cos \theta}{\sin \theta} b_\varphi
\end{aligned}
$$ 
and, therefore, the first equation of the system \eqref{BVK.19} becomes
$$
b_\varphi \, a_{u\theta\theta} + b_{\varphi\theta} \, a_{u\theta}
- a_{u\varphi} \, b_{\theta\theta} - b_\theta \, a_{u\varphi\theta} 
= \frac{\cos\theta}{\sin\theta} (a_{u\theta} \, b_\varphi - b_\theta \, a_{u\varphi}). 
$$
This is equivalent to saying 
$$
\big( a_{u\theta} \, b_\varphi - b_\theta \, a_{u\varphi} \big)_\theta
= 
\frac{\cos\theta}{\sin\theta}\big( a_{u\theta} \, b_\varphi - b_\theta \, a_{u\varphi} \big),
$$
which, after integration, leads to 
$$
a_{u\theta} \, b_\varphi - b_\theta \, a_{u\varphi} =  \sin \theta \, e^{G(\varphi)}
$$
for some function $G=G(\varphi)$.

As for the second equation in \eqref{BVK.19}, we have
$$
\begin{aligned}
& a_{u\theta} \, b_{\varphi\varphi} + b_\varphi \, a_{u\theta\varphi}
 - b_\theta \, a_{u\varphi\varphi} - a_{u\varphi} \, b_{\varphi\theta}
+ \sin \theta \cos \theta \, a_{u\theta} b_\theta
-  \sin \theta \cos \theta \, b_\theta a_{u\theta} 
 \\
  & = \frac{\cos \theta}{\sin \theta}\big( a_{u\varphi} \, b_\varphi -
 b_\varphi \, a_{u\varphi} \big),
 \end{aligned}
$$
so that
$$ 
\big( a_{u\theta} \, b_\varphi - b_\theta \, a_{u\varphi} \big)_\varphi = 0
$$
and $ a_{u\theta} \, b_\varphi - b_\theta \, a_{u\varphi}$ is a function of 
$\theta$ alone. Therefore, $G=G(\varphi)$ is a constant and we find
$$
a_{u\theta} \, b_\varphi - b_\theta \, a_{u\varphi} =  C \, \sin \theta
$$
for some scalar $C =C(u)$.

The condition above can be written in terms of the original components of $f_u$ and $X$,
 yielding
$$
f_u^\theta X^\varphi - f_u^\varphi X^\theta = \frac{C}{\sin \theta}
$$
and, for the embedded components,
$$
\tilde f_u^2 \tilde X^1 - \tilde f_u^1 \tilde X^2 = C.
$$
Since $X$ and $f_u$ are both tangent to the sphere,
the last equation means that there exists a constant $C$ such that
$$
X \times f_u = C \, \nn
$$
for all $x \in \Sbold^2$, so 
$$
| X \times f_u | = C, \quad x \in \Sbold^2.
$$
This formula describes vector fields that satisfy the constraint but may not
be globally smooth.  

Imposing now that $X$ is smooth, we deduce that it must vanish at least once on the sphere,
as must the flux $f_u$. By continuity, we necessarily have $C = 0$ so that 
$f_u$ and $X$ are colinear. This means that there exists a function $C(u, x)$
such that $f_u = C \, X$ for all $x \in \Sbold^2$. 
Finally, substituting $f_u$ for $C \, X$ in  $\Lie_X f_u=0$ yields $X(C) =0$.
This completes the proof of Lemma~\ref{BVK-8}.
\end{proof}

\begin{proof}[Proof of Theorem~\ref{TVS-1}] The first claim of the theorem is an immediate consequence
of Lemma~\ref{BVK-8}.

The equations are written along any integral curve of the vector field $X$.
To see this, fix a point $x \in \Sbold^2$ and introduce a parametrization $s \mapsto \varphi(s)$ of the (unique, integral)
curve passing through $x$ and 
such that the velocity vector $\frac{d}{ds}\varphi(s)$ equals the vector $X$ for every $s$.
Then define the function 
$$
\tilde u(t, s) := u(t, \varphi(s)), \qquad t\ge 0, \quad s\in \RR.
$$
Note that since $X(C) =0$ we have $C(\tilde u(t, s), \varphi(s)) = :\widetilde C(\tilde u)$ for all $s$. Therefore, using \eqref{TVS.2} and \eqref{BV.3},
we deduce that 
$$
\begin{aligned}
\del_s \widetilde C(\tilde u(t, s)) &= \widetilde C_u(\tilde u(t,s)) X(\varphi(s))(u)
	= du \big(\widetilde C_u(u(t, \varphi(s))) X(\varphi(s)) \big)
\\
	&= \nablac \big( f( u(t, \varphi(s)), \varphi(s)) \big)
	= -\del_t u(t, \varphi(s)) = - \del_t \tilde u (t, s),
\end{aligned}
$$
thus 
$$
\del_t u(t, s) + \del_s \widetilde C(u(t, s)) = 0.
$$
Note that this is a \emph{one-dimensional} conservation law in which 
the flux $\tilde C$ does not explicitly depend on the spatial coordinate. By Theorem~\ref{BV-3}, these conservation laws
admit total variation bounds of the form \eqref{BV.2}, provided that $f$ is divergence-free.
\end{proof}


\section{Finite volume method}   
\label{FVM-0}

The rest of this paper is devoted to establishing the convergence of the finite volume
method on a $d$-dimensional Riemannian manifold $(\MM,g)$.
We will construct approximate solutions to the Cauchy problem associated with the 
conservation law \eqref{IN.1} and the initial condition ($u_0 \in L^1(\MM) \cap L^\infty(\MM)$)  
\be
u(0,x) = u_0(x), \quad x \in \MM, 
\label{FVM.1a}
\ee
using a finite volume method based on monotone numerical flux-functions, and 
we will establish that such approximations converge strongly to the unique entropy solution of the problem.
We need not assume that the flux $f$ is divergence free (see \eqref{IN.2}), 
and on the map $\nablac f$
we require the following \emph{growth condition} 
\be
\label{FVM.4.1}
|\nablac f (u, x)| \le C_1 + C_2 |u|, \qquad (u,x) \in \RR \times M, 
\ee
for some constants $C_1,C_2>0$. 

Let us introduce a triangulation $\TT^h$ of $\MM$, that is, a set of 
curved polyhedra $K \subset \MM$ whose vertices are joined by geodesic lines.
We assume that $\MM = \bigcup_{K\in\TT^h} \overline K$ and that for every distinct elements $K_1,K_2\in \TT^h$, 
the set $K_1\cap K_2$ is either a common face of $K_1, K_2$ or else a submanifold with 
dimension at most $(d-2)$. To each element $K$ we associate the set of its faces $e$ by $\dK$, and 
we denote by $K_e$ the unique element distinct from $K$ sharing the face $e$ with $K$. 
The outward unit normal to an element $K$ at some point $x \in e$ is denoted by $\nek(x) \in T_xM$. 
Finally, $|K|$ and $|e|$ denote the $d$- and $(d-1)$-dimensional Hausdorff measures of $K$ and $e$, 
respectively, and $p_K = \sum_{e\in\dK}|e|$ is the measure of the boundary of $K$. 
When $N$ is a submanifold of $M$, we denote by $dv_N$
the volume element on $N$ induced by the metric $g$.

We impose that 
\be
h:=\sup_{\substack{x, y \in K \\ K\in \TT^h}} d_g(x,y) \to 0, 
\label{FVM.0.0}
\ee 
where $d_g(\cdot, \cdot)$ is the distance function associated to the metric $g$.
In particular, elements may become flat in the limit as $h$ tends to zero. 
Interestingly enough, this condition on $h$ is slightly more general than the corresponding condition in
\cite{CCL}, where the exterior diameter of the elements was considered.
We denote the time increment of the discretization by $\tau=\tau(h)$, and set $t_n:=\tau n$ for $n=0,1, \ldots$. 
We also assume that, as $h$ tends to zero,
\be
\label{FVM.0.1}
\tau \to 0, \qquad  h^2/\tau \to 0.
\ee

To each element $K$ and each edge $e\in \dK$ we associate a family of
locally Lipschitz continuous, numerical flux-functions $\fek : \RR\times \RR \to \RR$ satisfying the following properties.
\begin{itemize}
\item{\emph{Consistency property :} 
	\be
	\label{FVM.1}
	\fek(u,u) = \frac{1}{|e|} \int_e g\big(f(u,x), \nek(x)\big)\, dv_e, \qquad u\in \RR. 
	\ee
	}
\item{\emph{Conservation property :} 
	\be
	\label{FVM.2}
	\fek(u,v) = -\feke(v,u), \qquad u,v \in \RR. 
	\ee
	}
\item{\emph{Monotonicity property :}
	\be
	\label{FVM.3}
	\frac{\del}{\del u} \fek \ge 0, \quad \frac{\del}{\del v} \fek \le 0.
	\ee
	}
\end{itemize}

Note that, in the Euclidian case and when the flux depends only on the conservative variable,
 the consistency condition \eqref{FVM.1} reads $\fek(\ubar,\ubar)
 = f(\ubar)\cdot \nek$, which is the usual consistency condition (c.f.~\cite{CCL}).
For the sake of stability we impose the following \emph{CFL condition}, for every element $K\in \TT$:
\be
\label{FVM.5}
\begin{split} 
& \frac{\tau p_K}{|K|} \max_{K \in \TT^h, e \in \del K} \sup_{u,v \in \RR, u \neq v} \Big|{f_{e,K}(u,v) - f_{e,K}(v,v) \over u-v} \Big| 
< 1,
\\
& \frac{\tau p_K}{|K|} \, \Lip(f) <1, 
\end{split}
\ee 
where 
$$
\Lip(f) := \sup_{\substack {x\in M \\ \ubar \in \RR}} \Big| \frac{df}{du}(\ubar, x) \Big |_g 
$$
and $|\cdot|_g$ represents the norm induced on each tangent space $T_x M$ by the metric $g$.

As in the Euclidian case, the finite volume method can be motivated by formally averaging the 
conservation law \eqref{IN.1} over an element $K$ and using Stokes' theorem:
\be 
\begin{split} 
0 & = \frac{d}{dt} \frac{1}{|K|}\int_K u \, dv_M + 
\frac{1}{|K|}\int_K \nablac f(u,x)\, dv_M 
\\
& = \frac{d}{dt} \frac{1}{|K|} \int_K u \, dv_M  
+ \frac{1}{|K|}\int_{\dK} g\big( f(u,x), \nn(x) \big)\, dv_{\dK}
\\
& = \frac{d}{dt} \frac{1}{|K|} \int_K u \, dv_M + \frac{1}{|K|}  \sum_{e\in \dK}
	 \int_e g\big( f(u,x), \nek(x) \big)\, dv_e.
\end{split}
\nonumber
\ee 
We then introduce the approximation 
$$
\unk \approx \frac{1}{|K|} \int_K u \, dv_M, \quad 
\fek(\unk, \unke) \approx \frac{1}{|e|} \int_e  g\big( f(u,x), \nek(x) \big)\, dv_e
$$
and discretize the time derivative with a two-point scheme. 

Based on the above the finite volume approximations are defined by 
\be
\label{FVM.4}
\unnk := \unk - \sum_{e\in \dK} \frac{\tau |e|}{|K|} \fek
	(\unk, \unke).
\ee
The initial condition \eqref{FVM.1a} gives us 
\be 
u^0_K := \frac{1}{|K|} \int_K u_0(x)\, dv_M. 
\label{FVM.4.a}
\ee
Finally we define the function $u^h : \RR_+ \times M \to \RR$ by
\be 
u^h(t, x) = \unk, \quad x\in K, \quad t_n \le t \le t_{n+1}.
\label{FVM.4.b}
\ee

We can now state our main convergence result.

\begin{theorem}
\label{FVM-1}
Consider the Cauchy problem \eqref{IN.1}-\eqref{FVM.1a} associated with a conservation law on the Riemannian manifold
$(M,g)$ and some initial data $u_0 \in L^1(\MM)\cap L^\infty(\MM)$. 
Suppose that the flux satisfied the growth condition \eqref{FVM.4.1}. 
Let $\TT^h$ be a triangulation and $\tau=\tau(h)$ be a time increment satisfying the conditions \eqref{FVM.0.0}-\eqref{FVM.0.1}, 
and let $\fek$ be a family of numerical flux-functions satisfying the consistency, conservation, and monotonicity 
conditions \eqref{FVM.1}--\eqref{FVM.5}.
Let $u^h:\RR_+\times \MM \to \RR$ be the sequence of approximate solutions constructed by the finite volume method 
\eqref{FVM.4}--\eqref{FVM.4.b}. 
Then, for all $T>0$, the sequence $u^h$ is uniformly bounded in $L^\infty([0, T), L^1(\MM)\cap L^\infty(\MM))$ 
and converges almost everywhere towards the unique entropy solution $u \in L^\infty(\RR_+, L^1(\MM)\cap L^\infty(\MM))$
of the Cauchy problem \eqref{IN.1}-\eqref{FVM.1a}. 
\end{theorem}

The proof of this theorem is the subject of the next section.

\begin{remark} 
Without loss of generality, we may assume that $\Lip(f) >0$.
Indeed, in the trivial case $\Lip(f) = 0$,  $f$ does not depend on $u$ and the 
conservation law \eqref{IN.1} has the explicit solution 
$u(t, x) = u_0(x) - t \nablac f(x)$. Thus, for \emph{any} triangulation,
if we set $\fek := \frac{1}{|e|} \int_e g( f(x), \nn(x) )\, dv_e$, we immediately 
find from \eqref{FVM.4} that $\unk = u^0_K - n\tau \nablac f(x)$ for all $K$. This gives $u^h(t, x) =
u^h(0, x) - t \nablac f(x)$, and so $\lim_{h\to 0} u^h(t,x) = u_0(x) - t \nablac f(x)$,
 which is  the exact solution.
\end{remark}

\section{Proof of convergence}   
\label{PC-0}

This section contains the proof of Theorem~\ref{FVM-1}. 
We follow the strategy of proof proposed by Cockburn, Coquel, and LeFloch in the Euclidian case \cite{CCL}. 
As we will see, the derivation of the required estimates for a conservation law on a manifold is significantly more involved. 
 
\subsection{Discrete maximum principle, entropy dissipation and entropy inequalities} 

We begin with the discrete maximum principle and the $L^1$ contraction principle. 

\begin{lemma}
\label{PC-1}
The finite volume approximations satisfy 
\be
\label{PC.1}
\max_{K\in\TT^h} |\unk| \le
	 \big( \max_{K\in\TT^h} |u^0_K| + C_1 t_n \big)\, e^{C_2 t_n}, 
\ee
and, given another approximation sequence $v^h$ (defined in the obvious way from a different initial data
$v_0 \in L^1(\MM)\cap L^\infty(\MM)$), 
\be
\label{PC.2}
\| u^h (t + \tau) - v^h (t + \tau) \|_{L^1(M)}
	\le \| u^h (t) - v^h (t) \|_{L^1(M)}.
\ee
\end{lemma}
\begin{proof}
Let us first point out the following discrete version of Gronwall's inequality:  
if $w_n$ is a non-negative sequence satisfying $w_0 \le a$ and
\be
\label{PC.3}
w_n \le a + b \sum_{j=1}^{n-1} w_j, \qquad n \geq 1, 
\ee
for some non-negative constants $a, b$, then it follows that 
\be
\label{PC.4}
w_n \le a \, e^{bn}.
\ee
Namely, let us define $\phi_n$ by $\phi_0 = a$ and
$$
\phi_n := a + b \sum_{j=1}^{n-1} \phi_j.
$$
We have $\phi_n - \phi_{n-1} = b \phi_{n-1}$,
or $\phi_n = (1 + b)\phi_{n-1}.$ By induction we get $\phi_n = a (1 + b)^n$, and 
the desired estimate now follows from $1 + b \le e^b$ and $w_n \le \phi_n$.

Let us turn to the proof of \eqref{PC.1}. Using \eqref{FVM.4} and the consistency
condition \eqref{FVM.1}, we easily write
\be
\label{PC.5}
\begin{aligned}
& \unnk  = \bigg( 1 + \frac{\tau}{|K|} \sum_{e \in \dK} 
	\frac{ \fek(\unk, \unke) - \fek(\unk, \unk)}{\unke - \unk} |e| \bigg)\unk
\\
&  - \frac{\tau}{|K|} \sum_{e \in \dK}  \, \frac{\fek(\unk, \unke) 	
	-\fek(\unk, \unk)}{\unke - \unk} \, |e| \, \unke
	- \frac{\tau}{|K|} \int_K \nablac f(\unk, x) \, dv_M.
\\
\end{aligned}
\ee
The monotonicity condition \eqref{FVM.3} implies
$$
- \frac{\tau}{|K|} \sum_{e \in \dK} 
	\frac{ \fek(\unk, \unke) - \fek(\unk, \unk)}{\unke - \unk} \ge 0,
$$
so that, using the CFL condition \eqref{FVM.5},
$$
- \frac{\tau |e|}{|K|} \sum_{e \in \dK} 
	 \frac{ \fek(\unk, \unke) - \fek(\unk, \unk)}{\unke - \unk} 
\le \frac{\tau \, p_K}{|K|} \Lip(f) \le 1.
$$
Thus, we see that $\unnk$ is a convex combination of $\unk$ and
$(\unke)_{e\in\dK}$, up to a divergence term in \eqref{PC.5}. Therefore, we can
write
$$
\begin{aligned}
& \unnk \le \max (\unk, \max_{e\in \dK} \unke) 
	- \frac{\tau}{|K|} \int_K \nablac f(\unk, x) \, dv_M,
\\
& \unnk \ge \min (\unk, \min_{e\in \dK} \unke) 
	- \frac{\tau}{|K|} \int_K \nablac f(\unk, x) \, dv_M,
\end{aligned}
$$
and, by induction,
$$
\begin{aligned}
 |\unk| &\le \max_{K \in \TT^h} |u^0_K|
	+ \sum_{j=0}^{n-1} \tau \max_{K \in \TT^h}
		 \frac {1}{|K|} \int_K |\nablac f(u^j_K, x)| \, dv_M
\\
& \le \max_{K \in \TT^h} |u^0_K| +  \sum_{j=0}^{n-1} \tau \max_{K \in \TT^h}
	\big( C_1 + C_2 |u_K^j| \big) 
\\
& \le \max_{K \in \TT^h} |u^0_K| + n \tau C_1 + 
	\sum_{j=0}^{n-1} \tau C_2 \max_{K \in \TT^h} |u_K^j|
\end{aligned}
$$
for all $n> 0$. Using Gronwall's inequality \eqref{PC.3}-\eqref{PC.4} we arrive at the desired inequality \eqref{PC.1}.
 
We now turn to the proof of the $L^1$ contraction property \eqref{PC.2}. It will follow immediately once 
we check that 
\be
\label{PC.6}
\sum_{K \in \TT^h} |\unnk - v_K^{n+1}| |K|
 \le \sum_{K \in \TT^h} |\unk - v_K^{n}| |K|.
\ee
Introduce the discrete solution operator $H$,
$$
u^{n+1}\equiv(\unnk)_{K\in\TT^h}=(H_K(u^n))_{K\in\TT^h}=:H(u^n).
$$
Using the CFL and monotonicity conditions \eqref{FVM.5} and \eqref{FVM.3},
it is easy to check that $H$ is a monotone operator:
$$
\big( \unk \le  v^n_K, \quad K \in \TT^h \big)
\Rightarrow \big( H_K(\unk) \le H_K(v^n_K), \quad K \in \TT^h \big).
$$
Setting $a \wedge b = \max(a, b)$ and $a \vee b = \min(a,b)$, we have
$H_K( \unk \vee v_K^n) = \unnk \vee v_K^{n+1}$ and
$H_K( \unk \wedge v_K^n) = \unnk \wedge v_K^{n+1}$. 
Next, using the definition \eqref{FVM.4} we find
$$
\begin{aligned}
\unnk \vee v_K^{n+1} & -\unnk \wedge v_K^{n+1}
	= \unk \vee v_K^n - \unk \wedge v_K^n
\\
-\sum_{e\in\dK} & \frac{\tau |e|}{|K|} \big( \fek(\unk \vee 
	v_K^n, (\unk \vee v_K^n)_e) - \fek(\unk \wedge v_K^n,
		(\unk \wedge v_K^n)_e) \big).
\end{aligned}
$$

Observe now that
$$
(\unk \vee v_K^n)_e \le \unke \vee v^n_{K_e} \text{ and }
(\unk \wedge v_K^n)_e \ge \unke \wedge v^n_{K_e}
$$
This, together with $\del \fek/ \del v \le 0$, implies that
$$
\begin{aligned}
& \fek(\unk \vee 
	v_K^n, (\unk \vee v_K^n)_e) - \fek(\unk \wedge v_K^n,
		(\unk \wedge v_K^n)_e) 
\\
& \qquad \ge \fek(\unk \vee 
	v_K^n, \unke \vee v_{K_e}^n) - \fek(\unk \wedge v_K^n,
		\unke \wedge v_{K_e}^n),
\end{aligned}
$$
so that 
$$
\begin{aligned}
\unnk \vee v_K^{n+1} & -\unnk \wedge v_K^{n+1}
	\le \unk \vee v_K^n - \unk \wedge v_K^n
\\
-\sum_{e\in\dK} & \frac{\tau |e|}{|K|} \big( \fek(\unk \vee 
	v_K^n, \unke \vee v_{K_e}^n) - \fek(\unk \wedge v_K^n,
		\unke \wedge v_{K_e}^n) \big).
\end{aligned}
$$

The inequality \eqref{PC.6} now follows by multiplying by $|K|$ and
summing in $K \in \TT^h$, since the conservation property \eqref{FVM.2} allows us to cancel 
all of the numerical flux terms. The proof of Lemma \ref{PC.1} is complete.
\end{proof}

The discrete entropy inequalities will be formulated in terms of the following 
convex decomposition of $\unk$.
For each $K$, $e$, define (cf. \cite{CCL}, formula (3.6))
\be
\label{PC.7}
\tunnkee := \unk - \frac{\tau p_K}{|K|} \big( \fek(\unk, \unke)-
\fek(\unk, \unk) \big) 
\ee
and set 
\be
\label{PC.8}
\unnkee := \tunnkee - \frac{\tau}{|K|} \int_K \nablac 
	f( \unk, x) \, dv_M.
\ee
In view of \eqref{FVM.4} and the consistency property \eqref{FVM.1}, we have 
\be
\label{PC.9}
\unnk = \frac{1}{p_K} \sum_{e\in\dK} |e| \unnkee.
\ee

The following lemma establishes the existence of numerical flux terms  
and a first entropy inequality relating $\tunnkee$ and $\unk$. 
 We omit the proof and refer to \cite{CCL} and the references therein.

\begin{lemma}
Let $(U,F)=(U(u),F(u,x))$ be a convex entropy pair. Then there exists a family of numerical
flux function $\Fek : \RR^2 \to \RR$ satisfying the following conditions: 
\begin{itemize}
	\item{$\Fek$ is consistent with the entropy flux $F$: 
		$$
		\Fek (u, u) =
			\frac{1}{|e|} \int_e g\big( F(u,x), \nek(x) \big)\, dv_e, \qquad u \in \RR. 
		$$}
	\item{Conservation property : 
		$$
		\Fek(u, v) = - \Feke (v,u), \qquad u,v \in \RR. 
		$$}
	\item{Discrete entropy inequality:
		\be
		\label{PC.10}
		U(\tunnkee) - U(\unk) + \frac{\tau p_K}{|K|} \Big( \Fek(\unk, \unke)
			- \Fek(\unk, \unk) \Big) \le 0.
		\ee
		}
\end{itemize}
\end{lemma}

We can now write the discrete entropy inequality in terms of $\unnkee$ and $\unk$
(see~\eqref{PC.10}) :
\be
\label{PC.11}
\begin{split} 
& U(\unnkee) - U(\unk) + \frac{\tau p_K}{|K|} \Big( \Fek(\unk, \unke)
			- \Fek(\unk, \unk) \Big) \le \rnnkee,
\\
& \rnnkee := U(\unnkee) - U(\tunnkee).
\end{split}
\ee

The following entropy dissipation estimate is the key to the proof of convergence, 
as was observed in the Euclidian case by Cockburn, Coquel, and LeFloch \cite{CCL}. 

\begin{proposition}
\label{PC-3}
Let $U:\RR \to \RR$ be a strictly convex, smooth function. Then, for all $n$ we have
\be
\label{PC.12}
\begin{aligned}
& \sum_{K\in \TT^h} U(\unnk) |K| + \frac{\alpha}{2} 
	\sum_{\substack{K \in \TT^h \\ e \in \dK}} 
	\frac{|e| |K|}{ p_K } \big| \unnkee - \unnk \big|^2
\\
& \quad \le \sumke U(\unk) |K| +  \sumke \tau |e| \Fek( \unk, \unk)
	+\sumke \frac{|e| |K|}{p_K} \rnnkee, 
\end{aligned}
\ee
where $\alpha$ denotes the modulus of convexity of $U$.
\end{proposition}

\begin{proof}
First, we observe that from an elementary lemma on convex functions (Lemma 3.5 of
\cite{CCL})
and \eqref{PC.9} one easily obtains
\be
\label{PC.13}
 \sumk |K| U(\unnk) + \frac{\alpha}{2} \sumke \frac{|e| |K|}{p_K} 
	\big| \unnkee - \unnk \big|^2 
 \le \sumke \frac{|e| |K|}{p_K} U(\unnkee).
\ee
Next, we sum up the inequality \eqref{PC.11} with respect to all $K$ and $e$, and divide by $p_K$ to obtain
\be
\label{PC.14}
 \sumke \frac{|e| |K|}{p_K} U(\unnkee) - \sumk U(\unk) |K|  - 
 	\sumke \tau |e| \Fek (\unk, \unk) - \sumke \frac{|e| |K|}{p_K} \rnnkee \le 0. 
\ee
Note that, again, all of the flux terms $\Fek(\unk, \unke)$ cancel out due to conservation.
Combining \eqref{PC.14} with \eqref{PC.13} gives \eqref{PC.12}. The proof of 
Proposition \ref{PC-3} is complete.
\end{proof}

As a consequence of the entropy dissipation estimate \eqref{PC.12}, one can deduce 
the following bound: 
\be
\label{PC.15}
\sum_{n=1}^N \sumke \frac{|e| |K|}{p_K} \big| u^n_{K,e} - \unk \big|^2
	\le \|u_0\|^2_{L^2(M)} + C(T),
\ee
where $C(T)$ is a constant depending only on $T$ and $N$ is an integer
such that $T -N \tau < \tau$.
For instance, it suffices to take $U(u) = u^2/2$ in \eqref{PC.12}.

We now derive a local entropy inequality
which does not involve the numerical flux  and is a consequence of the inequality \eqref{PC.11}. 
The proof is completely similar to that of Lemma~3.6 in \cite{CCL} and we omit it.
This inequality will then allow us to prove a global discrete entropy inequality. 

\begin{lemma}
Let $(U,F)$ be a convex entropy pair. Then for all elements $K$ and all faces 
$e \in \dK$, we have
\be
\label{PC.16}
\begin{aligned}
& \frac{|K|}{p_K} U(\unnkee) - \frac{|K|}{p_K} U(\unk)
	+\frac{|K_e|}{p_{K_e}} U(\unnkeee) - \frac{|K_e|}{p_{K_e}} U(\unke)
\\
&\quad +\frac{\tau}{|e|} \int_e g \big( F(\unke, x) - F(\unk, x), \nek (x) \big)\, dv_e
	\le \frac{|K|}{p_K} \rnnkee + \frac{|K_e|}{p_{K_e}} R^{n+1}_{K_e,e},
\end{aligned}
\ee
where $\rnnkee$ is defined in \eqref{PC.11}.
\end{lemma}

\begin{proposition}
\label{PC-5}
Let $(U,F)$ be a convex entropy pair and let $\phi = \phi(t, x) 
\in C_c \big( [0, T) \times \MM \big)$ be a test function.
 For each element $K$ and each face
 $e$ of $K$, set
\be
\pne = \frac{1}{\tau |e|} \int_{t_n}^{t_{n+1}} \!\!\!		
	\int_e \phi(t,x) \,dv_e dt, \qquad \phnk = \sum_{e \in \dK} 
	\frac{|e|}{p_K} \pne,
\label{PC.17}
\ee
\be 
\widehat{\del_t \phi_K}^n = \frac{1}{\tau} (\phnk - \hat \phi_K^{n-1}).
\label{PC.18}
\ee
Then, $\unk$ satisfies the following inequality:
\be
\begin{aligned}
- \sum_{n=0}^N &\sumk \int_{t_n}^{t_{n+1}} 
	\int_K \Big( U(\unk) \widehat{\del _t \phi_K}^{n} +
	g \big( F(\unk, x), \nabla_g \phi(t,x) \big) \Big) \,dv_M dt
\\
& + \sumk \int_K U(u_K^0) \hat \phi_K^0 \,dv_M
	\le E^h +  Q^h + R^h, 
\end{aligned}
\label{PC.19}
\ee
where the error terms $E^h$, $Q^h$ and $R^h$ are given by
\be
\label{PC.20}
E^h = \sum_{n=0}^N \sumke
	\frac{|e| |K|}{p_K} U(\unnkee) (\phnk - \pne),
\ee
\be
\label{PC.21}
\begin{aligned}
Q^h &:  = \sum_{n=0}^N \sumk Q^n_K,
\\
Q^n_K & :=\sum_{e\in \dK}
		\tau \int_e g\big( F(\unk, x), \nek(x) \big) \, dv_e\,
		\pne - \int_{t_n}^{t_{n+1}} \int_K g\big( F(\unk, x), 
	\nabla_g \phi(t,x) \big)\, dv_K dt,
\end{aligned}
\ee
\be
\label{PC.22}
\begin{split} 
& R^h = \sum_{n=0}^N \sumke \frac{|K| |e|}{p_K} \pne \rnnkee 
\\
&  R_{K,e}^{n+1} = U(\unnkee) - U(\tunnkee).
\end{split}
\ee
\end{proposition}

\begin{proof}
Multiply the inequality \eqref{PC.16} by $\pne |e|$ and sum over all $K \in \TT^h$, 
$e \in \dK$. Note that the first four terms in \eqref{PC.16}, as well as the 
last, give the same result whether they involve $K$ or $K_e$
(since the sum is taken over the same set of elements). We thus obtain
$$
\begin{aligned}
& 2 \sumke \frac{|K| |e|}{p_K} \big( U(\unnkee) - U(\unk) \big) \pne
\\
& \quad + \sumke \tau |e| \pne \Big( \frac{1}{|e|} \int_e g \big( F(\unke, x) - F(\unk, x)
	, \nek(x) \big)\, dv_e \Big) \le 2 \sumke \frac{|K| |e|}{p_K} \rnnkee \pne.
\end{aligned} 
$$
Also, due to the conservation and consistency properties of $\Fek$ we have
$$
\sumke \frac{1}{|e|} \int_e g \big(F(\unk, x), \nek (x) \big) \, dv_e 
= - \sumke \frac{1}{|e|} \int_e g \big( F(\unke, x), \nek (x) \big) \, dv_e,
$$
which implies
\be
\label{PC.23}
\begin{aligned}
& \sumke \frac{|K| |e|}{p_K} \big( U(\unnkee) - U(\unk) \big) \pne
\\
& \quad - \sumke \tau |e| \pne \frac{1}{|e|} \int_e g \big( F(\unk, x)
	, \nek(x) \big)\, dv_e \le  \sumke \frac{|K| |e|}{p_K} \rnnkee \pne.
\end{aligned}
\ee
Observe now that, setting $Q_K^n$ as in \eqref{PC.21}, we deduce 
from the definition of $\pne$ (see~\eqref{PC.17}) that 
\be
\label{PC.24}
\begin{aligned}
& \sumke \tau |e| \pne \frac{1}{|e|} \int_e g \big( F(\unk, x)
	, \nek(x) \big)\, dv_e
\\
& = \sumk \int_{t_n}^{t_{n+1}}\int_K g \big( F(\unk, x), \nabla \phi(t,x) \big)\, dv_M dt
	+ \sumk Q_K^n.
\end{aligned}
\ee
From the convexity of $U$ and Jensen's inequality, we see that
\be
\label{PC.25}
\sumk |K| U( \unnk) \phnk \le
	\sumke \frac{|K| |e|}{p_K} U( \unnkee) \phnk.
\ee
Thus, from \eqref{PC.23}, \eqref{PC.24} and \eqref{PC.25} and after summation in $n$,
we obtain 
$$
\begin{aligned}
& \sum_{n=0}^N \sumk |K| \big( U(\unnk) - U(\unk) \big) \phnk -
	\sum_{n=0}^N \sumk  \int_{t_n}^{t_{n+1}} \int_K
	g \big( F(\unk, x), \nabla \phi(t,x) \big) \, dv_M dt
\\
& \le E^h + \sum_{n=0}^N \sumk Q_K^n + 
	\sum_{n=0}^N \sumke \frac{|K| |e|}{p_K} \pne \rnnkee,
\end{aligned}
$$
with $E^h$ given by \eqref{PC.20}. Summing by parts in the first term, one obtains
the desired entropy inequality \eqref{PC.19}. The proof of Proposition~\ref{PC-5} is complete.
\end{proof}


\subsection{Measure-valued solutions and convergence analysis}

We now prove the strong convergence of the approximate solutions
$u^h$ given by the finite volume method towards the entropy solution to the Cauchy problem under consideration. 
To this end, we need to rely on the framework of measure-valued solutions to conservation laws developed in \cite{DiPerna} 
in the Euclidian case and extended to manifolds by 
Ben-Artzi and LeFloch \cite{BenArtziLeFloch}.
The key question is to justify the passage to the limit as $h \to 0$ in the 
inequality \eqref{PC.19} and use the entropy dissipation estimate \eqref{PC.15} to control the error terms.

It is well-known that to the sequence $u^h$ (which is uniformly bounded in 
$L^\infty( [0,T) \times M)$ for every $T>0$) we can associate a subsequence and a
Young measure $\nu : [0, T) \times M \to \mathrm{Prob}(\RR)$,
which is a family of probability measures in $\RR$ parametrized by $(t,x) \in [0, T) \times M$. 
The Young measure allows us to determine all weak-$*$ limits of composite functions $a(u^h)$, 
for arbitrary real continuous functions $a$, according to the following property :
\be
\label{PC.26}
a(u^h) \mathrel{\mathop{\rightharpoonup}\limits^{*}} \langle \nu, a \rangle
\quad \text{ when } h \to 0
\ee
where we use the notation 
$$
\langle \nu, a \rangle := \int_\RR a(\lambda) \, d\nu(\lambda).
$$

In view of the above property, the passage to the limit in the 
left hand side of \eqref{PC.19} is (almost) immediate. The key uniqueness theorem 
\cite{DiPerna, BenArtziLeFloch} tells us that once we know that $\nu$ is a 
measure-valued solution to the conservation law, we can prove that the support
of each probability measure $\nu_{t,x}$ actually reduces to a single point 
$u(t, x)$, if the same is true for $t=0$, that is, $\nu_{t,x}$ is the Dirac measure $\delta_{u(t,x)}$. 
It is then standard to deduce that the convergence in \eqref{PC.26} is actually strong, 
and that, in particular, $u^h$ converges strongly to $u$ which in turn is the unique 
entropy solution of the Cauchy problem under consideration.

\begin{lemma}
\label{PC-6}
Let $\nu : [0, T) \times M \to \mathrm{Prob}(\RR)$ be a Young measure 
associated with the sequence $u^h$. 
Then, for every convex entropy / entropy-flux pair $(U,F)$ and for every 
non-negative, smooth function $\phi : [0,T) \times M \to \RR_+$ with compact support,
one has 
\be
\label{PC.27}
\begin{aligned}
& - \int_0^T \int_M \langle \nu_{t,x}, U(\cdot) \rangle \del_t \phi(t,x) +
	g\big( \langle \nu_{t,x}, F(\cdot, x) \rangle, \nabla \phi(t,x) \big)\, dv_M dt +
	\int_M U(u_0(x)) \phi(0, x) \, dv_M
\\
& \quad \le  \int_0^T \int_M \langle \nu_{t,x}, \nablac F (\cdot, x) \rangle
	 \phi(t,x) dv_M dt
	- \int_0^T \int_M \langle \nu_{t,x}, U'(\cdot) \nablac f( \cdot, x)
	 \rangle \phi(t, x) \, dv_M dt.
\end{aligned}
\ee
\end{lemma}

\begin{proof}
First, we claim that the discretizations \eqref{PC.17}, \eqref{PC.18} 
satisfy the following estimates, whose proofs are elementary:
$$
\sup_{\substack{t \in [t_n, t_{n+1}] \\ x \in K}} \big| \del_t\phi (t,x) 
	- \widehat{\del_t \phi_k}^n \big| \le (\tau +h_K) \| \phi \|_{	
	C^2([t_n, t_{n+1}] \times K)},
$$
and, for each face $e$,
$$
|\pne - \phnk| \le (\tau +h_K) \| \phi \|_{
	C^1([t_n, t_{n+1}] \times K)}.
$$
In view of these estimates and \eqref{PC.26}, the left hand side of \eqref{PC.19}
converges to the left-hand side of \eqref{PC.27}. 

We now consider the error terms on the right-hand side of \eqref{PC.19}. First, note 
that $E^h$ equals the corresponding one in \cite{CCL}, so we 
refer the reader to that paper for a proof that $E^h$ converges to zero. Note 
that the use of estimate \eqref{PC.15} is essential in that proof.

Next, let us prove that 
\be
\label{PC.28}
\lim_{h\to 0} R^h = - \int_0^T \int_M \langle \nu_{t,x}, U'(\cdot)\,
	\nablac f (\cdot, x) \rangle \phi(t, x)\, dv_M dt. 
\ee
Recall that $R^h$ is given by \eqref{PC.22}. Introducing the averages
$$
\phi_K^n := \frac{1}{|K| \tau} \int_{t_n}^{t_{n+1}}
	\int_K \phi(t, x)\, dv_M dt,
$$
we claim that it suffices to prove that 
\be
\lim_{h\to 0} \Big( R^h + \sum_n \sumk  U'(\unk)
	\tau \int_K \nablac f (\unk, x)\, dv_M \, \pnk \Big) = 0.
\label{PC.29}
\ee
Indeed, we have
$$
\begin{aligned}
& - \sum_n \sumk  U'(\unk)
	{\tau} \int_K \nablac f (\unk, x)\, dv_M \, \pnk
\\
& \qquad = - \sum_n \sumk U'(\unk) \int_{t_n}^{t_{n+1}} \int_K
	\nablac f(\unk, x) \phi(t, x)\, dv_M dt
\\
& \qquad \qquad - \sum_n \sumk U'(\unk) \Big( \inttn \int_K
	\nablac f(\unk, x)\, dv_M \frac{1}{|K|} \int_K \phi(t, x)\, dv_M dt
\\
& \hspace{45mm} - \inttn \int_K \nablac f(\unk, x) \phi(t,x)\, dv_M dt\Big),
\end{aligned}
$$
so that taking the limit as $h \to 0$ we find, for the first term,
$$
\begin{aligned}
& \lim_{h\to 0} - \int_0^T \int_M U'(u^h) \nablac f(u^h, x) \phi(t,x)\, dv_M dt
\\
& \qquad = - \int_0^T \int_M \langle \nu_{t,x}, U'(\cdot)\,
	\nablac f (\cdot, x) \rangle \phi(t, x)\, dv_M dt,
\end{aligned}
$$
while the second term easily converges to zero. 

We now prove \eqref{PC.29}.
First, observe that from the definition of $\pnk$, it is easy to see that we may 
as well consider the term $R^h$ with $\pnk$ instead of $\pne$. We thus have
$$
\begin{aligned}
 R^h  &+ \sum_{n=0}^N \sumke \frac{|K| |e|}{p_K}  U'(\unk)
	\frac{\tau}{|K|} \int_K \nablac f (\unk, x)\, dv_M \, \pnk 
\\
& = R^h + \sum_{n=0}^N \sumke \frac{|K| |e|}{p_K}  U'(\unnk)
	\frac{\tau}{|K|} \int_K \nablac f (\unnk, x)\, dv_M \, \phi_K^{n+1}
\\
& \qquad + \sumke \frac{|K| |e|}{p_K}  U'(u^0_K)
	\frac{\tau}{|K|} \int_K \nablac f (u^0_K, x)\, dv_M \, \phi_K^{0}
\\
& = \sum_{n=0}^N \sumke \frac{|K| |e|}{p_K} \pnk \Big( 
	U(\unnkee) - U(\tunnkee) + U'(\unnk)\frac{\tau}{|K|} \int_K
	\nablac f(\unnk, x)\, dv_M \Big)
\\
& \qquad + \sum_{n=0}^N \sumke \frac{|K| |e|}{p_K}  U'(\unnk)
	\frac{\tau}{|K|} \int_K \nablac f(\unnk, x)\, dv_M  (\phi_K^{n+1}	
	- \pnk)
\\
& \qquad \qquad - \sumk U'(u_K^0) \tau \int_K \nablac f(u_K^0, x)\, dv_M \phi_K^0
 = A + B + C. 
\end{aligned}
$$
Note that all sums in $n$ may be written from $0$ to $N$ because of the compact
support of $\phi$. 

Consider the first term, $A$. From the bound 
\eqref{FVM.4.1} on the divergence of $f$ and from the maximum principle 
\eqref{PC.1} we conclude that
$$
\frac{\tau}{|K|} \int_K \nablac f (\unnk, x)\, dv_M = \Ocal(\tau).
$$ 
Thus, Taylor developing $U$ around the point $\tunnkee$, and using \eqref{PC.8}, 
the definition of $\unnkee$, we find 
$$
\begin{aligned}
|A| & \le \sum_n \sumke \frac{|K| |e|}{p_K} \Big| U'(\unnkee) \Ocal(\tau) 
	+ \Ocal(\tau^2) - U'(\unnk) \Ocal(\tau) \Big| \pnk
\\
& \le \Ocal(\tau) \sup U'' \sum_n \sumke \frac{|K| |e|}{p_K}|\unnkee - \unnk|\pnk.
\end{aligned}
$$
Next, using the estimate \eqref{PC.15} and the Cauchy-Schwartz inequality, we find
that the term inside the sum above is bounded, so that $A = \Ocal(\tau)$.

Consider now the second term, $B$. The function $\phi$ is regular, so that
$\phi_K^{n+1} - \pnk = \Ocal(\tau)$. Therefore, we have
$$
|B| \le \sum_n \sumk |K| \Ocal(\tau^2) = \Ocal(\tau).
$$
Finally, using again the bound \eqref{FVM.4.1}, we see that the last term $C$
 tends to zero with $\tau$. This concludes the proof of \eqref{PC.28}.

Next, consider the term $Q^h=\sum_{n,K} Q_K^n$
given by \eqref{PC.21}. Using Green's formula, we can write
$$
\begin{aligned}
Q_K^n & = \sum_{e\in\dK} \frac{1}{|e|}
	\int_e g\big( F(\unk, x), \nek(x) \big) \, dv_e
		\int_{t_n}^{t_{n+1}} \int_e \phi(t,x) \, dv_e dt
\\
& - \int_{t_n}^{t_{n+1}} \int_{\dK} g\big( F(\unk, x) \phi(t, x),
	\nn(x) \big)\, dv_{\dK} dt 
	+ \int_{t_n}^{t_{n+1}} \int_K  \nablac F (\unk, x) \phi(t,x)\, dv_M dt.
\end{aligned}
$$
We wish to obtain
\be
\label{PC.30}
\lim_{h\to 0} Q^h =
 \int_0^T \int_M \langle \nu_{t,x},  \nablac F (\cdot, x) \rangle \phi(t, x)
	\, dv_M dt,
\ee
so that according to the above equality we must prove that
\be
\label{PC.31}
\begin{aligned}
&\lim_{h\to 0} \sum_n \sumke \Big( \frac{1}{|e|}
	\int_e g\big( F(\unk, x), \nek(x) \big) \, dv_e
		\int_{t_n}^{t_{n+1}} \int_e \phi(t,x) \, dv_e dt
\\
& \qquad - \int_{t_n}^{t_{n+1}} \int_{e} g\big( F(\unk, x) \phi(t, x),
	\nek(x) \big)\, dv_{\dK} dt \Big) = 0.
\end{aligned}
\ee
Let $x_e$ be any point on the edge $e$. Note that
$$
\begin{aligned}
 \sum_n & \sumke  
	\inttn \int_e g\big( F(\unk, x), \nek(x) \big) \, 
		\Big( \frac{1}{|e|} \int_e \phi(t,x) \, dv_e\, - \phi(t, x) \Big) dv_e dt
\\
& = \sum_n \sumke  
	\inttn \int_e \Big( g\big( F(\unk, x), \nek(x) \big) 
	- g \big( F(\unk, x_e), \nek(x_e)\big) \Big) \cdot
\\ 
& \hspace{65mm} \cdot \Big( \frac{1}{|e|} \int_e \phi(t,x) \, dv_e\, 
	- \phi(t, x) \Big) dv_e dt,
\end{aligned}
$$
so that taking absolute values and in view of the regularity
of $F$ and $\phi$, we can bound the last term by
$$
C T \sumk p_K h^2.
$$
using the CFL condition \eqref{FVM.5}, we 
can finally write
$$
\begin{aligned}
\Big| Q^h  - \int_0^T \int_M   \nablac F  (u^h, x) \phi(t,x) \,dv_M dt \Big|
	&\le \frac{C T}{\Lip(f)} \sumk p_K h^2 
\\
	&\le C \sumk \frac{|K|}{\tau} h^2 
		= \frac{h^2}{\tau} C |M|,
\end{aligned}
$$
so that taking the limit as $h\to 0$ and using \eqref{FVM.0.1} yields
$$
Q^h \to \int_0^T \int_M \langle \nu_{t,x},  \nablac F (\cdot, x) \rangle \phi(t, x)
	\, dv_M dt,
$$
which is \eqref{PC.30}.
This completes the proof of Lemma \ref{PC-6}.
\end{proof}

\begin{remark} From the proof of the previous lemma, we see that if 
the divergence-free condition \eqref{IN.2} is satisfied, then the error
term $Q^h$ converges to zero with $h$ and $R^h$ vanishes altogether.
In that case, we obtain, instead of \eqref{PC.27}, the usual weak formulation
$$
 - \int_0^T \int_M \langle \nu_{t,x}, U(\cdot) \rangle \del_t \phi(t,x) +
	g\big( \langle \nu_{t,x}, F(\cdot, x) \rangle, \nabla \phi(t,x) \big)\, dv_M dt +
	\int_M U(u_0(x)) \phi(0, x) \, dv_M \le 0.
$$
\end{remark}

\begin{proof}[Proof of Theorem \ref{FVM-1}]
According to the inequality \eqref{PC.27}, we have for all convex entropy pairs
$(U,F)$,
$$
 \del_t \langle \nu, U(\cdot) \rangle + \nablac (\langle \nu, F(\cdot, x)\rangle)
 \le \langle \nu, (\nablac F)(\cdot, x)\rangle - \langle \nu, U'(\cdot)
	(\nablac f)(\cdot, x)\rangle
$$
in the sense of distributions in $[0, T)\times M$.
Since, for $t = 0$, the Young measure $\nu$ is the Dirac mass $\delta_{u_0}$ 
(because $u_0$ is a bounded function), we know from \cite{BenArtziLeFloch}
 that there exists a unique function $u \in L^\infty([0, T)\times M)$ such that the measure $\nu$ remains the 
Dirac mass $\delta_u$ for all time $0 \le t \le T$. Moreover,
this implies that the approximations $u^h$ converge strongly to $u$
on compact sets at least. This concludes the proof.
\end{proof}


\section*{Acknowledgments} 
P.A. was supported by the Portuguese Foundation for Science 
and Technology (FCT) through grant SFRH/BD/17271/2004.
P.G.L. was partially supported by the Centre National de la Recherche Scientifique (CNRS)
and by a grant (Number 2601-2) from the Indo-French Center for the Promotion of Advanced Research (IFCPAR).  


\newcommand{\auth}{\textsc}


\begin{thebibliography}{10}  

\bibitem{BenArtziLeFloch} \auth{Ben Artzi M. and LeFloch P.G.,}
{\sl Hyperbolic conservation laws on manifolds. The well-posedness theory,}
in preparation. 
 
\bibitem{CCL} \auth{Cockburn B., Coquel F., and LeFloch P.G.,}
{\sl Convergence of finite volume methods for multidimensional conservation laws,}
SIAM J. Numer. Anal. {\bf 32} (1995), 687--705.
  
\bibitem{CCLS} \auth{Cockburn B., Coquel F., LeFloch P.G., and Shu C.W.,}
{\sl Convergence of finite volume methods for multidimensional conservation laws,}
Preprint, Institute for Mathematics and its Applications, Minneapolis, 1989. 
  
\bibitem{DiPerna} \auth{DiPerna R.J.,}
{\sl Measure-valued solutions to conservation laws,}  
Arch. Rational Mech. Anal. {\bf 88} (1985), 223--270.  

\bibitem{Kruzkov} \auth{Kruzkov S.,}
{\sl First-order quasilinear equations with several space variables} (in Russian), Mat. USSR Sb. {\bf 123} (1970), 228--255;  
English Transl. in Math. USSR Sb. {\bf 10,} 217--243. 


\end{thebibliography}
\end{document}